%% file: 2019-03-03-from_mkc.tex
\documentclass[journal]{IEEEtran}

\ifCLASSINFOpdf
\else
\fi
\usepackage{amsfonts}
\usepackage{mathtools}
\usepackage{subcaption}
\usepackage{amssymb}
\usepackage{bbm}
\usepackage{mathtools}
\usepackage{enumitem}
\usepackage{cite}

\usepackage{tikz}
\usepackage{tkz-graph}
\usepackage{graphicx}      

\usepackage[american]{circuitikz}
\makeatletter
\ctikzset{bipoles/vsources/minus/sign/rotation/.initial=0}
\pgfcircdeclarebipole{}{\ctikzvalof{bipoles/vsourceam/height}}{vsourceAM}{\ctikzvalof{bipoles/vsourceam/height}}{\ctikzvalof{bipoles/vsourceam/width}}{
	\pgfsetlinewidth{\pgfkeysvalueof{/tikz/circuitikz/bipoles/thickness}\pgfstartlinewidth}
	\pgfpathellipse{\pgfpointorigin}{\pgfpoint{0}{\pgf@circ@res@up}}{\pgfpoint{\pgf@circ@res@left}{0}}
	\pgfusepath{draw}
	\ifpgf@circ@oldvoltagedirection
	\pgftext[bottom,rotate=90,y=\ctikzvalof{bipoles/vsourceam/margin}\pgf@circ@res@down]{$+$}
	\pgftext[top,rotate=90,y=\ctikzvalof{bipoles/vsourceam/margin}\pgf@circ@res@up]{$-$}
	\else
	\pgftext[bottom,rotate=90,y=\ctikzvalof{bipoles/vsourceam/margin}\pgf@circ@res@down]{\rotatebox{\ctikzvalof{bipoles/vsources/minus/sign/rotation}}{$-$}}
	\pgftext[top,rotate=90,y=\ctikzvalof{bipoles/vsourceam/margin}\pgf@circ@res@up]{$+$}
	\fi
}

\pgfcircdeclarebipole{}{\ctikzvalof{bipoles/cvsourceam/height}}{cvsourceAM}{\ctikzvalof{bipoles/cvsourceam/height}}{\ctikzvalof{bipoles/cvsourceam/width}}{
	\pgfsetlinewidth{\pgfkeysvalueof{/tikz/circuitikz/bipoles/thickness}\pgfstartlinewidth}
	\pgfpathmoveto{\pgfpoint{\pgf@circ@res@left}{\pgf@circ@res@zero}}
	\pgfpathlineto{\pgfpoint{\pgf@circ@res@zero}{\pgf@circ@res@up}}
	\pgfpathlineto{\pgfpoint{\pgf@circ@res@right}{\pgf@circ@res@zero}}
	\pgfpathlineto{\pgfpoint{\pgf@circ@res@zero}{\pgf@circ@res@down}}
	\pgfpathlineto{\pgfpoint{\pgf@circ@res@left}{\pgf@circ@res@zero}}
	\pgfusepath{draw}
	
	\ifpgf@circ@oldvoltagedirection
	\pgftext[bottom,rotate=90,y=\ctikzvalof{bipoles/cvsourceam/margin}\pgf@circ@res@left]{$+$}
	\pgftext[top,rotate=90,y=\ctikzvalof{bipoles/cvsourceam/margin}\pgf@circ@res@right]{$-$}
	\else
	\pgftext[bottom,rotate=90,y=\ctikzvalof{bipoles/cvsourceam/margin}\pgf@circ@res@left]{\rotatebox{\ctikzvalof{bipoles/vsources/minus/sign/rotation}}{$-$}}
	\pgftext[top,rotate=90,y=\ctikzvalof{bipoles/cvsourceam/margin}\pgf@circ@res@right]{$+$}
	\fi
}
\makeatother

\input{ka-newcommands}

\renewcommand{\EP}{\hfill\IEEEQED}

\newtheorem{thm}{Theorem}
\newtheorem{defn}[thm]{Definition}
\newtheorem{lem}[thm]{Lemma}
\newtheorem{lemma}[thm]{Lemma}
\newtheorem{cor}[thm]{Corollary}
\newtheorem{ex}[thm]{Example}

\newtheorem{remark}[thm]{Remark}

\newtheorem{problem}[thm]{Problem}

\usepackage{tcolorbox}

\usepackage[colorlinks=true,linkcolor=green]{hyperref}
\hypersetup{
	colorlinks,
	linkcolor={blue!80!black},
	citecolor={green!80!black},
	urlcolor={brown!80!black}
}
\def\equationautorefname~#1\null{Equation~(#1)\null}
\usetikzlibrary{calc,patterns,decorations.pathmorphing,decorations.markings}

\newcounter{todocounter}
\usepackage[colorinlistoftodos]{todonotes}

\makeatletter

\newcommand{\Rmnum}[1]{\expandafter\@slowromancap\romannumeral #1@}
\makeatother

\usepackage{pmat}

\begin{document}
%
\title{A Unifying Framework for Strong Structural Controllability}

\author{\IEEEauthorblockN{Jiajia Jia,
Henk J. van Waarde,
Harry L. Trentelman, 
and
M. Kanat Camlibel
}

\thanks{The authors are with the Bernoulli Institute for Mathematics, Computer Science and Artificial Intelligence, University of Groningen, The Netherlands (e-mail: {\footnotesize{\tt j.jia@rug.nl},{\tt h.j.van.waarde@rug.nl},
	{\tt h.l.trentelman@rug.nl}, {\tt m.k.camlibel@rug.nl})}}
%

}


%



\IEEEtitleabstractindextext{%
%
%
%
%
%

\begin{abstract}
This paper deals with strong structural controlla\-bility of linear systems. In contrast to existing work, the structured systems studied in this paper have  a so-called zero/nonzero/arbitrary structure, which means that some of the entries are equal to zero, some of the entries are arbitrary but nonzero, and the remaining entries are arbitrary (zero or nonzero). We formalize this in terms of pattern matrices whose entries are either fixed zero, arbitrary nonzero, or arbitrary. 
We establish necessary and sufficient algebraic conditions for strong structural controllability in terms of full rank tests of certain pattern matrices. We also give a necessary and sufficient graph theoretic condition for the full rank property of a given pattern matrix. This graph theoretic condition makes use of a new color change rule that is introduced in this paper. 
Based on these two results, we then establish a necessary and sufficient graph theoretic condition for strong structural controllability. Moreover, we relate our results to those that exists in the literature, and explain how our results generalize previous work.
\end{abstract}

\begin{IEEEkeywords}
Strong structural controllability, Network controllability, Structured system, Pattern matrices.
\end{IEEEkeywords}}

\maketitle

\IEEEdisplaynontitleabstractindextext

%
\IEEEpeerreviewmaketitle

\section{Introduction}
Controllability is a fundamental concept in systems and control. For linear time-invariant systems of the form 
\begin{equation}
\label{linearsystem}
\dot{x} = Ax + Bu,
\end{equation}
controllability can be be verified using the  Kalman rank test or the Hautus test \cite{TSH2012}. 
Often, the exact values of the entries in the matrices $A$ and $B$ are not known, but the underlying interconnection structure between the input and state variables is known exactly.  

In order to formalize this, Mayeda and Yamada have introduced a framework in which, instead of a fixed pair of real matrices, only the \emph{zero/nonzero} structure of $A$ and $B$ is given \cite{MY1979}. This means that each entry of these matrices is known to be either a \emph{fixed zero} or an \emph{arbitrary nonzero} real number. Given this zero/nonzero structure, they then study controllability of the family of systems for which the state and input matrices have this zero/nonzero structure. In this setup, this family of systems is called \emph{strongly structurally controllable} if all members of the family are controllable in the classical sense \cite{MY1979}.

To the best of our knowledge, all existing literature up to now (except for \cite{MZC2014}) has considered strong structural controllability under the above basic assumption that for each of the entries of the system matrices there are only {\em two possibilities}: it is either a fixed zero, or an arbitrary nonzero value \cite{MY1979,Bachmann1981,RSW1992,OMP1993,JSA2011,CM2013,RHS2014,TD2015}.

There are, however, many scenarios in which, in addition to these two possibilities, there is a third possibility, namely, that a given entry is not a fixed zero or nonzero, but can take {\em any real value}.  In such a scenario, it is not possible to represent the system using a zero/nonzero structure, but a third possibility needs to be taken into account.  To illustrate this, consider the following example.

\begin{figure}[h!]
	\centering
	\scalebox{0.8}{
	\begin{circuitikz}[scale=0.8]
		\draw
		(3,-1)--(-1,-1);
		\draw (-1,2) to[american voltage source, V<= $V$] (-1,-1);
		\draw (-1,2)
		to [R = $R$, *-, i<=$I_{R}$] (3,2)
		to [C = $C_1$, -*,v=$V_{C_1}$,i= $I_{C_1}$] (3,-1)--(7,-1)
		to [american current source, l_ = $I$] (7,2);
		\draw (3,2) to [L = $L$,*-*,i=$I_{\textsc{L}}$] (7,2)--(7,5.25);
		\draw (-1,2) -- (-1,5.25) to [C = $C_2$,v=$V_{C_2}$] (3,5.25);
		\draw (7,2) -- (7,5.25);
		\draw (3,5.25)
		to [american controlled voltage source, l = $GI_{C_1}$,bipoles/vsources/minus/sign/rotation=90 ] (7,5.25) ;
	\end{circuitikz}}
	\caption{Example of electrical circuit.}
	\label{g:electriccircuit}
\end{figure}
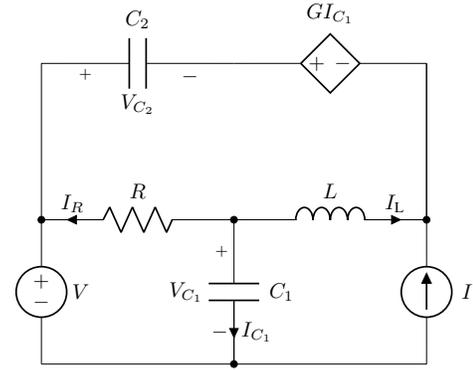
\begin{ex}\label{ex:electriccircuit}
	The electrical circuit in Figure \ref{g:electriccircuit}
	consists of a resistor, two capacitors, an inductor, an independent voltage source, an independent current source and a current controlled voltage source.  Assume that the parameters $R, C_1, C_2$ and $L$ are positive but not known exactly. We denote the current through $R$, $L$, and $C_1$ by $I_{R}$, $I_L$, and $I_{C_1}$, respectively, and the voltage across $C_1$ and $C_2$ by $V_{C_1}$ and $V_{C_2}$, respectively. 
	The current controlled voltage source is represented by $G I_{C_1}$ with gain $G$ assumed to be positive. 
	Define the state vector as $x = [V_{C_1} ~ V_{C_2} ~ I_L]^T$ and the input as $u = [V ~ I]^T$. By Kirchhoff's current and voltage laws, the circuit is represented by a linear time-invariant system \eqref{linearsystem} with 
		\begin{equation}\label{e:e1}
	A = \begin{bmatrix}
-\frac{1}{R C_1} & 0 & -\frac{1}{C_1}\\[1mm]
0 & 0 & -\frac{1}{C_2}\\[1mm]
\frac{R-G}{RL} & \frac{1}{L} & - \frac{G}{L}
	\end{bmatrix}, ~~ B =  
	\begin{bmatrix}
	\frac{1}{R C_1} & 0\\[1mm]
	0 & -\frac{1}{C_2} \\[1mm] 
	 \frac{G-R}{RL} & 0
	\end{bmatrix}.
	\end{equation}
	
	Recall that the parameters $R,C_1,C_2,L > 0$ are not known exactly. This means that the matrices in \eqref{e:e1} are not known exactly, but we do know that they have the following structure. Firstly, some entries  are \emph{fixed zeros}. Secondly, some of the entries are always \emph{nonzero}, for instance, the entry with value $-\frac{1}{R C_1}$. The third type of entries, those with value $\frac{R-G}{RL}$ and $\frac{G-R}{RL}$, can be either zero (if $R = G$) or nonzero.
	Since the system matrices in this example do not have a zero/nonzero structure, the existing tests for strong structural controllability \cite{MY1979,Bachmann1981,RSW1992,OMP1993,JSA2011,CM2013,RHS2014,TD2015} are not applicable. 
\end{ex}

A similar problem as in Example \ref{ex:electriccircuit} appears in the context of linear networked systems. Strong structural controllability of such systems has been well-studied \cite{CM2013, TD2015,MZC2014,MHM2018,PP2013}. In the setup of these references, the weights on the edges of the network graph are unknown, while the network graph itself is known. Under the assumption that the edge weights are arbitrary but nonzero, linear networked systems can thus be regarded as systems with a given zero/nonzero structure. This zero/nonzero structure is determined by the network graph, in the sense that nonzero entries in the system matrices correspond to edges in the network graph. However, often even exact knowledge of the network graph is not available, in the sense that it is unknown whether certain edges in the graph exist or not.  This issue of missing knowledge appears, for example, in social networks \cite{K2006}, the world wide web \cite{RAJ2015}, biological networks \cite{CMN2008,RM2009} and ecological systems \cite{TE2012}. Another cause for uncertainty about the network graph might be malicious attacks and unintentional failures. This issue is encountered in transportation networks \cite{VM2005}, sensor networks \cite{KM2007} and gas networks\cite{CB2009}. 

\begin{ex}\label{ex:network}
Consider a network of three agents with single-integrator dynamics, represented by
	\[\dot{x}_i = v_{i}\]
	for $i = 1,2,3$. Here $x_i \in \mathbb{R}$ is the state of agent $i$ and $v_i \in \mathbb{R}$ is its input. The communication between the agents is represented by the graph in Figure \ref{g:network}.
	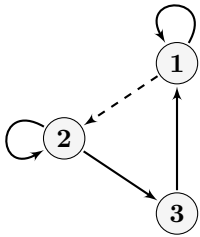
\begin{figure}[h!]
		\centering
		\centering
		\begin{tikzpicture}[scale=0.5]
		\tikzset{VertexStyle1/.style = {shape = circle,
				color=black,
				fill=white!96!black,
				minimum size=0.5cm,
				text = black,
				inner sep = 2pt,
				outer sep = 1pt,
				minimum size = 0.55cm}
		}	
		\tikzset{VertexStyle2/.style = {shape = circle,
				ball color = black!80!yellow,
				text = white,
				inner sep = 2pt,
				outer sep = 0pt,
				minimum size = 10 pt}}
		\node[VertexStyle1,draw](1) at (0,2) {$\bf 1$};
		\node[VertexStyle1,draw](2) at (-3,0) {$\bf 2$};
		\node[VertexStyle1,draw](3) at (0,-2) {$\bf 3$};
		%
		\Edge[ style = {->,> = latex',pos = 0.2},color=black
		, labelstyle={inner sep=0pt}](2)(3);
		\Edge[ style = {->,> = latex',pos = 0.7, dashed},color=black, labelstyle={inner sep=0pt}](1)(2);
		\Loop[style={> = latex',->, out=150, in=-150,line width=0.8pt,color=black}, dist=1.5cm](2)
		\Edge[style = {->,> = latex',pos = 0.7},color=black
		, labelstyle={inner sep=0pt}](3)(1);	
		\Loop[style={> = latex',->, out=60, in=120,line width=0.8pt,color=black}, dist=1.5cm](1);	
		\end{tikzpicture}
		\caption{Example of a networked system.}
		\label{g:network}
	\end{figure}
	The links $(1,1)$, $(2,2)$, $(2,3)$ and $(3,1)$ are known to exist, while the link $(1,2)$ is uncertain in the sense that it may or may not be present. This is represented by solid and dashed edges, respectively. Agents $1$ and $2$ are only affected by the states of their neighbors, while agent $3$ is also influenced by an external input $u \in \mathbb{R}$. This means that $v_1 = w_{11}x_1 + w_{13}x_3$, $v_2 = w_{21}x_1 + w_{22}x_2$ and $v_3 = w_{32}x_2 + u$. Here the weights $w_{11}, w_{22}, w_{32}$ and $w_{13}$ are \emph{nonzero} since they correspond to existing edges, while the weight $w_{21}$ that corresponds to the uncertain link is arbitrary (zero or nonzero).
	We can write the network system in compact form \eqref{e:e1} by defining 
	\begin{equation}
	\label{e:networkexample}	
	A = 	\bbm
	w_{11} & 0 & w_{13}\\
	w_{21} & w_{22} & 0\\
	0 & w_{32} & 0   \\
	\ebm , B = \bbm 0\\ 0\\ 1\ebm.
	\end{equation}
	Since $w_{21}$ can be zero or nonzero, the system matrices in \eqref{e:networkexample} do not have a zero/nonzero structure. 
\end{ex}

To conclude, both in the context of modeling physical systems, as well as in representing networked systems, capturing the system simply by a zero/nonzero structure is not always possible, and a more general concept of system structure is required. Therefore, in this paper we will extend the notion of zero/nonzero structure, and study strong structural controllability for families of systems having such more general structure. In particular, our main contributions are the following:
\begin{enumerate}
	\item We extend the notion of zero/nonzero structure to a more general {\em zero/nonzero/arbitrary structure}, and formalize this structure in terms of suitable pattern matrices. 
	\item We establish necessary and sufficient conditions for strong structural controllability for families of systems with a given zero/nonzero/arbitrary structure. These conditions are of an algebraic nature and can be verified by a rank test on two pattern matrices.
	\item We provide a graph theoretic condition for a given pattern matrix to have full row rank. This condition can be verified using a new \emph{color change rule}, that will be defined in this paper. 
	\item We establish a graph theoretic test for strong structural controllability for the new families of structured systems.
	\item Finally, we relate our results to those existing in the literature by showing how existing results can be recovered from those we present in this paper. We find that seemingly incomparable results of \cite{TD2015} and \cite{MZC2014} follow from our main results, which reveals an overarching theory. For these reasons, our paper can be seen as a unifying approach to strong structural controllability of linear time-invariant systems. 
\end{enumerate}

We conclude this section by giving a brief account of research lines that are related to strong structural controllability but that will not be pursued in this paper. The concept of weak structural controllability was introduced by Lin in \cite{Lin74} and has been studied extensively, see \cite{Lin74,SP1976,DCW2003,LSB2011,SH2013,CvdWB2017,CvdW2018}. Another, more recent, line of work focuses on structural controllability of systems for which there are \emph{dependencies} among the arbitrary entries of the system matrices \cite{LM2017,JTBC2018}. An important special case of dependencies among parameters arises when the state matrix is constrained to be symmetric, which was considered in \cite{MKBP2018,MBP2018,MHM2018}. The problem of \emph{minimal input selection} for controllability has also been well-studied, see, e.g., \cite{PKA2016,Olshevsky2014,SCL2016,TRPJ2016}. Finally, weak and strong structural \emph{targeted} controllability have been investigated in \cite{LCPPP2018} and \cite{MCT2015,vWCT2017}, respectively.

The outline of the rest of the paper is as follows. In Section~\ref{s:pre},  we present some preliminaries. Next, in Section~\ref{s:pro},  we formulate the main problem treated in this paper. Then, in Section~\ref{s:main} we state our main results. Section~\ref{s:dis} contains a comparison of our results with previous work. In Section~\ref{S:proof} we state proofs of the main results. Finally, in Section~\ref{s:con} we formulate our conclusions. 

\section{Preliminaries}\label{s:pre}
Let $\mathbb{R}$ and $\mathbb{C}$ denote the fields of real and complex numbers, respectively. 
The spaces of $n$-dimensional real and complex vectors are denoted by $\mathbb{R}^{n}$ and $\mathbb{C}^{n}$, respectively. 
Likewise, the space of $n \times m$ real  matrices is denoted by $\mathbb{R}^{n \times m}$.

Moreover, $I$ and $0$ will denote the identity and zero matrix of appropriate dimensions, respectively. 
	
In this paper, we will use so-called pattern matrices. By a pattern matrix we mean a matrix with entries in the set of symbols $\{0,\ast,?\}$. These symbols will be given a meaning in the sequel.

The set of all $p \times q$ pattern matrices will be denoted by $\{0,\ast,?\}^{p \times q}$.
For a given $p \times q$ pattern matrix $\mathcal{M}$, we define the \emph{pattern class} of $\mathcal{M}$ as
\begin{equation*}
\begin{aligned}
\mathcal{P}(\mathcal{M}) := \{M \in \mathbb{R}^{p \times q} \mid  & M_{ij} = 0 \text{ if } \mathcal{M}_{ij} = 0, \\& M_{ij} \neq 0 \text{ if } \mathcal{M}_{ij} = \ast \}.
\end{aligned}
\end{equation*}
This means that for a matrix $M \in \mathcal{P}(\mathcal{M})$, the entry $M_{ij}$ is either (i) \emph{zero} if $\mathcal{M}_{ij} = 0$, (ii) \emph{nonzero} if $\mathcal{M}_{ij} = \ast$, or (iii) \emph{arbitrary} (zero or nonzero) if $\mathcal{M}_{ij} = \: ?$.

\section{Problem formulation}\label{s:pro}
Let $\mathcal{A} \in \{0,\ast,?\}^{n \times n}$ and $\mathcal{B} \in \{0,\ast,?\}^{n \times m}$ be pattern matrices.
Consider the linear dynamical system
\begin{equation}\label{e:gss}		
		\dot{x}(t) = A x(t) + B u(t),
\end{equation}
where the system matrix $A$ is in $\mathcal{P}(\mathcal{A})$ and the input matrix $B$ is in $\mathcal{P}(\mathcal{B})$, and where $x \in \mathbb{R}^n$ is the state and $u \in \mathbb{R}^m$ is the input.

We will call the family of systems \eqref{e:gss} a {\em structured system}. To simplify the notation, we denote this structured system by the ordered pair of pattern matrices $(\mathcal{A},\mathcal{B})$.
\begin{ex}\label{ex:same circuit}
Consider the electrical circuit discussed in Example \ref{ex:electriccircuit}. Recall that this was modelled as the state space system \eqref{e:e1} in which the entries of the system matrix and input matrix were either fixed zeros, strictly nonzero or undetermined. This can be represented as a structured system $(\mathcal{A},\mathcal{B})$ with pattern matrices 
	\begin{equation} 
	\calA = \begin{bmatrix}
\ast & 0 & \ast\\
0 & 0 & \ast\\
? & \ast  & \ast
\end{bmatrix} \mbox{ and } \calB = 
\begin{bmatrix}
\ast & 0\\
0 & \ast \\ 
?& 0
\end{bmatrix}.
	\end{equation}
\end{ex}
In this paper we will study structural controllability of structured systems. In particular, we will focus on strong structural controllability, which is defined as follows.	
\begin{defn}\label{d:SSC}
The system $(\mathcal{A},\mathcal{B})$ is called \emph{strongly structurally controllable} if the pair $(A,B)$ is controllable for all $A \in \mathcal{P}(\mathcal{A})$ and $B \in \mathcal{P}(\mathcal{B})$. 
\end{defn}

The concept of strong structural controllability  was introduced by Mayeda and Yamada in the 1970's \cite{MY1979} and has been further investigated in \cite{RSW1992, JSA2011}. In these works, the structured system matrices $\mathcal{A}$ and $\mathcal{B}$ are restricted to only contain $0$ and $*$ entries. 
In the context of controllability of networked systems \cite{LSB2011}, the study of strong structural controllability was extended to linear networked  systems, see e.g., \cite{CM2013,MZC2014,TD2015}. 
In these references, a networked system is also represented by a linear structured system $(\mathcal{A},\mathcal{B})$ where $\mathcal{A}$ is determined by the structure of the network and $\mathcal{B}$ encodes the leader nodes through which external inputs are injected into the network. In this framework, a common assumption is that each input only affects a single node in the network. This means that $\mathcal{B}$ is a pattern matrix
with exactly one $*$ in each column and at most one $*$ in each row. In addition, the pattern matrices $\mathcal{A}$ studied in \cite{CM2013,MZC2014,TD2015} can be seen as special cases of the pattern matrices studied in the present paper. Indeed, the papers \cite{CM2013} and \cite{TD2015} consider the case in which $\mathcal{A}$ only contains $0$ and $*$ entries. Furthermore, the paper \cite{MZC2014} deals with pattern matrices $\mathcal{A}$ whose diagonal entries are all $?$ and none of the off-diagonal entries is $?$.

Up to now, a framework for studying strong structural controllability of $(\mathcal{A},\mathcal{B})$ where $\mathcal{A} \in \{0,*,?\}^{n \times n}$ and $\mathcal{B} \in \{0,*,?\}^{n \times m}$ are \emph{general} pattern matrices has not yet been developed. 
Therefore, the problem that we will investigate in the present paper is
stated as follows.
\begin{problem}
	Given two  pattern matrices $\mathcal{A} \in \{0,\ast,?\}^{n \times n}$ and $\mathcal{B} \in \{0,\ast,?\}^{n \times m}$, provide necessary and sufficient conditions under which $(\mathcal{A},\mathcal{B})$ is strongly structurally controllable.
\end{problem}

In the remainder of this paper, we will simply call the structured system $(\mathcal{A},\mathcal{B})$ {\em controllable} if it is strongly structurally controllable.

\begin{remark}\label{r:weak}
	In addition to strong structural controllability, in the past also \emph{weak structural controllability} has been studied extensively. This concept was introduced by Lin in \cite{Lin74}. Instead of requiring {\em all} systems in a family associated with a given structured system to be controllable, weak structural controllability only asks for the existence of at least one controllable member of that family, see \cite{Lin74,SP1976,DCW2003}. 
	In these references, conditions were established for weak structural controllability of structured systems in which the pattern matrices only contain $0$ or $?$ entries. The question then arises: is it possible to generalize the results from \cite{Lin74,SP1976,DCW2003} to structured systems in the context of our paper, with more general pattern matrices $\mathcal{A} \in \{0,\ast,?\}^{n \times n}$ and $\mathcal{B} \in \{0,\ast,?\}^{n \times m}$. 
	Indeed, it turns out that the results in \cite{Lin74,SP1976,DCW2003} can immediately be applied to assess weak structural controllability of our more general structured systems. To show this, for given pattern matrices $\mathcal{A} \in \{0,\ast,?\}^{n \times n}$ and $\mathcal{B} \in \{0,\ast,?\}^{n \times m}$ we define two new pattern matrices $\mathcal{A}'\in \{0,?\}^{n \times n}$ and $\mathcal{B}' \in \{0,?\}^{n \times m}$ as follows: $\mathcal{A}'_{ij} = 0 \iff \mathcal{A}_{ij} = 0$ and $\mathcal{B}'_{ij} = 0 \iff \mathcal{B}_{ij} = 0$. The new structured system $(\mathcal{A}',\mathcal{B}' )$  is now a structured system of the form studied in \cite{Lin74,SP1976,DCW2003}. Using the fact that weak structural controllability is a generic property \cite{SP1976}, it can then be shown that weak structural controllability of $(\mathcal{A}',\mathcal{B}' )$ is equivalent to that of $(\mathcal{A},\mathcal{B})$. In other words, weak structural controllability of general $(\mathcal{A},\mathcal{B})$ can be verified using the conditions established in previous work \cite{Lin74,SP1976,DCW2003}. 
\end{remark}

\section{Main results}\label{s:main}

In this section, the main results of this paper will be stated. Firstly, we will establish an algebraic condition for controllability of a given structured system. This condition states that controllability of a structured system is equivalent to full rank conditions on two pattern matrices associated with the system. Secondly, a graph theoretic condition for a given pattern matrix to have full row rank will be given in terms of a so-called \emph{color change rule}. Finally, based on the above algebraic condition and graph theoretic condition, we will establish a graph theoretic necessary and sufficient condition for controllability of a structured system.

Our first main result is a rank test for controllability of a structured system. In the sequel, we say that a pattern matrix $\calM$ {\em has full row rank\/} if every matrix $M\in\calP(\calM)$ has full row rank. 

\begin{thm}
	\label{t:algebraic} 
	The system $(\mathcal{A},\mathcal{B})$ is controllable if and only if the following two conditions hold:
	\begin{enumerate}[label = \arabic*.]
		\item The pattern matrix $\bbm \calA & \calB \ebm$ has full row rank.\\[-3mm]
		\item The pattern matrix $\bbm \bar{\calA} & \calB \ebm$ has full row rank where $\bar{\mathcal{A}}$ is the pattern matrix obtained from $\mathcal{A}$ by modifying the diagonal entries of $\mathcal{A}$ as follows:
		\begin{equation}\label{e:barofX}
		\bar{\mathcal{A}}_{ii} := \begin{cases}
		\ast & \text{if } \mathcal{A}_{ii} = 0,  \\
		? & \text{otherwise}.
		\end{cases}	
		\end{equation}
	\end{enumerate}
\end{thm}

We note here that the two rank conditions in Theorem \ref{t:algebraic} are independent, in the sense that one does not imply the other in general. To show that the first rank condition does not imply the second, consider the pattern matrices $\mathcal{A}$, the corresponding $\bar{\mathcal{A}}$, and $\mathcal{B}$ given by 
		\begin{equation*}
		\mathcal{A} = \begin{bmatrix}
		\ast   & \ast  \\
		0	   & 0\\
		\end{bmatrix}, \: \bar{\mathcal{A}} = \begin{bmatrix}
		?   & \ast  \\
		0	   & \ast\\
		\end{bmatrix} \mbox{ and }
		\mathcal{B} = \begin{bmatrix}
		\ast   \\
		\ast    \\
		\end{bmatrix}.
		\end{equation*}
		It is evident that the pattern matrix $\bbm \calA & \calB \ebm$ has full row rank.
		However, for the choice 
		\begin{equation*}
		\bar{A} = \bbm 0 & 1\\ 0 & 1 \\ \ebm \in \mathcal{P}(\bar{\mathcal{A}}) \mbox{ and } B = \bbm 1 \\ 1 \ebm \in \calP(\calB),
		\end{equation*}
		the matrix $\bbm \bar{A} & B \ebm$ does not have full row rank.	
		
		To show that the second condition does not imply the first one, consider the pattern matrix $\mathcal{A}$, the corresponding $\bar{\mathcal{A}}$, and $\mathcal{B}$ given by
		\[
		\mathcal{A} = \begin{bmatrix}
		?   & 0  \\
		\ast	   & 0\\
		\end{bmatrix}, \: \bar{\mathcal{A}} = \begin{bmatrix}
		?   & 0  \\
		\ast	   & \ast\\
		\end{bmatrix} \mbox{ and } \mathcal{B} = \begin{bmatrix}
		\ast   \\
		\ast    \\
		\end{bmatrix}.
		\]
Obviously, the pattern matrix $\bbm \bar{\calA} & \calB \ebm$ has full row rank. However, for the choice 
		\begin{equation*}
		A = \bbm 1 & 0\\ 1 & 0 \\ \ebm \in \mathcal{P}(\mathcal{A}) \mbox{ and } B = \bbm 1 \\ 1 \ebm \in \mathcal{P}(\mathcal{B}),
		\end{equation*}
		we see that $\bbm A & B \ebm$ does not have full row rank.

Next, we discuss a noteworthy special case in which the first rank condition in Theorem~\ref{t:algebraic} is implied by the second one. Indeed, if none of the diagonal entries of $\calA$ is zero, it follows from \eqref{e:barofX} that $\calP(\calA) \subseteq \calP(\bar{\calA})$. Hence, we obtain the following corollary to Theorem~\ref{t:algebraic}.
 
\begin{cor}\label{c:algebraic}
Suppose that none of the diagonal entries of $\calA$ is zero. Let $\bar{\calA}$ be as defined in \eqref{e:barofX}. The system $(\mathcal{A},\mathcal{B})$ is controllable if and only if  $\bbm \bar{\calA} & \calB \ebm$ has full row rank.
\end{cor}

Note that both $\bbm \mathcal{A} & \mathcal{B} \ebm$ and $\bbm\bar{\mathcal{A}} & \mathcal{B}\ebm$ appearing in Theorem~\ref{t:algebraic} are $n \times (n + m)$ pattern matrices. Next, we will develop a graph theoretic test for checking whether a given pattern matrix has full rank. 
To do so, we first need to introduce some terminology. 

Let $\mathcal{M} \in \{0,\ast,?\}^{p \times q}$ be a pattern matrix with $p \leq q$.  
We associate a directed graph $G(\mathcal{M}) = (V,E)$ with $\mathcal{M}$ as follows. Take as node set $V = \{1, 2, \ldots, q\}$ and define the edge set $E \subseteq V \times V$ such that $(j,i) \in E$ if and only if $\mathcal{M}_{ij} = \ast$ or $\mathcal{M}_{ij} =?$. If $(i,j) \in E$, then we call $j$ an {\em out-neighbor\/} of $i$. Also, in order to distinguish between $\ast$ and $?$ entries in $\mathcal{M}$, we define two subsets $E_\ast$ and $E_?$ of the edge set $E$ as follows: $(j,i) \in E_\ast$ if and only if $\mathcal{M}_{ij} = \ast$ and $(j,i) \in E_?$ if and only if $\mathcal{M}_{ij} = ?$. Then, obviously, $E = E_\ast \cup E_?$ and $E_\ast \cap E_? = \emptyset$.
To visualize this, we use  solid  and  dashed arrows to represent edges in $E_\ast$ and $E_?$, respectively.
\begin{ex}
As an example, consider the pattern matrix $\mathcal{M}$ given by 
\begin{equation*}
\mathcal{M} = \begin{bmatrix}
0	   & 0     & \ast  & 0  & 0 \\
0	   & \ast  & \ast  & ?  & \ast  \\
\ast   & 0     & ?     & 0  & 0  \\
0      & \ast  & 0     & 0  & ? \\
\end{bmatrix}.
\end{equation*}
The associated directed graph $G(\mathcal{M})$ is then given in Figure~\ref{g:illustrate}. 

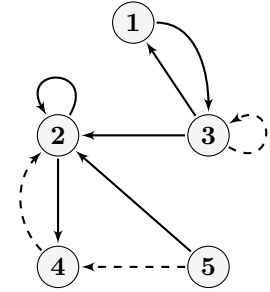
\begin{figure}[h!]
	\centering
	\centering
	\begin{tikzpicture}[scale=0.5]
		\tikzset{VertexStyle1/.style = {shape = circle,
	        color=black,
	        fill=white!96!black,
	        minimum size=0.5cm,
			text = black,
			inner sep = 2pt,
			outer sep = 1pt,
			minimum size = 0.55cm}
		}	
	\tikzset{VertexStyle2/.style = {shape = circle,
			ball color = black!80!yellow,
			text = white,
			inner sep = 2pt,
			outer sep = 0pt,
			minimum size = 10 pt}}
	\node[VertexStyle1,draw](1) at (0,3) {$\bf 1$};
	\node[VertexStyle1,draw](2) at (-2,0) {$\bf 2$};
	\node[VertexStyle1,draw](3) at (2,0) {$\bf 3$};
	\node[VertexStyle1,draw](4) at (-2,-3.5) {$\bf 4$};
	\node[VertexStyle1,draw](5) at (2,-3.5) {$\bf 5$};
	
	\Edge[ style = {->,> = latex',pos = 0.2,out=0,in=90},,color=black
	, labelstyle={inner sep=0pt}](1)(3);
		\Edge[ style = {->,> = latex',pos = 0.7},,color=black
	, labelstyle={inner sep=0pt}](3)(1);
	\Loop[style={> = latex',->, out=60, in=120,line width=0.8pt}, dist=1.5cm](2)
	\Edge[style = {->,> = latex',pos = 0.7},color=black
	, labelstyle={inner sep=0pt}](2)(4);

	\Edge[ style = {->,> = latex',pos = 0.3},,color=black
	, labelstyle={inner sep=0pt}](3)(2);
	\Loop[style={> = latex',->, out=-30, in=30,line width=0.8pt,dashed}, dist=1.5cm](3);		
	\Edge[ style = {->,> = latex',pos = 0.3,dashed, out = 135, in=-135},,color=black
	, labelstyle={inner sep=0pt}](4)(2);
	\Edge[ style = {->,> = latex',pos = 0.3,dashed},,color=black
	, labelstyle={inner sep=0pt}](5)(4);
	\Edge[ style = {->,> = latex',pos = 0.3},color=black
	, labelstyle={inner sep=0pt} ](5)(2);
	\end{tikzpicture}
	\caption{ The graph $G(\calM)$ associated with $\mathcal{M}$.}
	\label{g:illustrate}
\end{figure}
\end{ex}

Next, we introduce the notion of {\em colorability\/} for $G(\mathcal{M})$:
\ben
\item Initially, color all nodes of $G(\mathcal{M})$ white.
\item\label{i:color} If a node $i$ has exactly one white out-neighbor $j$ and $(i,j) \in E_*$, we change the color of $j$ to black.
\item Repeat step \ref{i:color} until no more color changes are possible$.$
\een
The graph $G(\mathcal{M}) $ is called {\em colorable\/} if the nodes $1, 2, \ldots, p$ are colored black following the procedure above. Note that the remaining nodes $p+1, \ldots ,q$ can never be colored black since they have no incoming edges. 

We refer to step \ref{i:color} in the above procedure as the {\em color change rule\/}. Similar color change rules have appeared in the literature before (see e.g. \cite{H2010,MZC2014,TD2015}). Unlike some of these rules, node $i$ in step \ref{i:color} does not need to be black in order to change the color of a neighboring node.



\begin{ex}			
	Consider the pattern matrix $\mathcal{M}$ given by
	\begin{equation*}
	\mathcal{M} = \begin{bmatrix}
	\ast	& 0   & 0    & 0    & \ast  & 0 \\
	0	    & ?   & 0    & \ast & 0     & \ast\\
	\ast	& 0   & 0    & \ast & 0     & 0 \\
	0	    & ?   & \ast & \ast & 0     & 0 \\
	\end{bmatrix}.
	\end{equation*}
	The directed graph $ G(\mathcal{M})$ associated with $\mathcal{M}$ is depicted in Figure \ref{g:G1}. By repeated application of the color change rule as shown in Figure \ref{g:G2} to \ref{g:G4}, we obtain the derived set $\mathcal{D} = \{1,2,3,4\}$. Hence, $G(\mathcal{M})$ is colorable.
	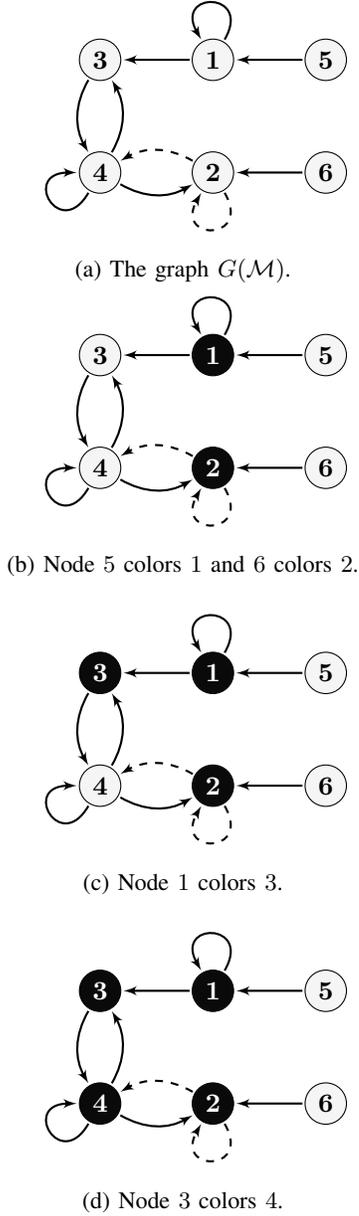
\begin{figure}[h!]
		\centering
		\begin{subfigure}{0.3\textwidth}
			\centering
			\begin{tikzpicture}[scale=0.5]
		\tikzset{VertexStyle1/.style = {shape = circle,
	        color=black,
	        fill=white!96!black,
	        minimum size=0.5cm,
			text = black,
			inner sep = 2pt,
			outer sep = 1pt,
			minimum size = 0.55cm}
		}		

			\node[VertexStyle1,draw](1) at (0,3) {$\bf 1$};
			\node[VertexStyle1,draw](2) at (0,0) {$\bf 2$};
			\node[VertexStyle1,draw](3) at (-3,3) {$\bf 3$};
			\node[VertexStyle1,draw](4) at (-3,0) {$\bf 4$};
			\node[VertexStyle1,draw](5) at (3,3) {$\bf 5$};
			\node[VertexStyle1,draw](6) at (3,0) {$\bf 6$};
			\Edge[ style = {->,> = latex',pos = 0.2},,color=black
			, labelstyle={inner sep=0pt}](5)(1);
			\Edge[style = {->,> = latex',pos = 0.7},color=black
			, labelstyle={inner sep=0pt}](1)(3);
			\Edge[style = {->,> = latex',pos = 0.7},color=black
			, labelstyle={inner sep=0pt}](6)(2);
			\Edge[ style = {->,> = latex',pos = 0.7,dashed,out=150,in=30},,color=black
			, labelstyle={inner sep=0pt}](2)(4);
			\Edge[ style = {->,> = latex',pos = 0.7,out=-30,in=-150},,color=black
			, labelstyle={inner sep=0pt}](4)(2);
			\Edge[ style = {->,> = latex',pos = 0.3,out=60,in=-60},,color=black
			, labelstyle={inner sep=0pt}](4)(3);
			\Edge[ style = {->,> = latex',pos = 0.3,out=-120,in=120},,color=black
			, labelstyle={inner sep=0pt}](3)(4);
			\Loop[style={> = latex',->, out=60, in=120,line width=0.8pt}, dist=1.5cm](1)
			\Loop[style={> = latex',->, out=-60, in=-120,line width=0.8pt,dashed}, dist=1.5cm](2)
			\Loop[style={> = latex',->, out=-120, in=180,line width=0.8pt}, dist=1.5cm](4)
			\end{tikzpicture}
			\caption{\centering The graph $G(\mathcal{M})$.}
			\label{g:G1}
		\end{subfigure}
		\vspace{.3cm}
		\begin{subfigure}{0.3\textwidth}
			\centering
			\begin{tikzpicture}[scale=0.5]
		\tikzset{VertexStyle1/.style = {shape = circle,
	        color=black,
	        fill=white!96!black,
	        minimum size=0.5cm,
			text = black,
			inner sep = 2pt,
			outer sep = 1pt,
			minimum size = 0.55cm}
		}		
			\tikzset{VertexStyle2/.style = {shape = circle,
	        color=black,
	        fill=black!96!white,
	        minimum size=0.5cm,
			text = white,
			inner sep = 2pt,
			outer sep = 1pt,
			minimum size = 0.55cm}
}
			\node[VertexStyle2,draw](1) at (0,3) {$\bf 1$};
			\node[VertexStyle2,draw](2) at (0,0) {$\bf 2$};
			\node[VertexStyle1,draw](3) at (-3,3) {$\bf 3$};
			\node[VertexStyle1,draw](4) at (-3,0) {$\bf 4$};
			\node[VertexStyle1,draw](5) at (3,3) {$\bf 5$};
			\node[VertexStyle1,draw](6) at (3,0) {$\bf 6$};
			\Edge[ style = {->,> = latex',pos = 0.2},,color=black
			, labelstyle={inner sep=0pt}](5)(1);
			\Edge[style = {->,> = latex',pos = 0.7},color=black
			, labelstyle={inner sep=0pt}](1)(3);
			\Edge[style = {->,> = latex',pos = 0.7},color=black
			, labelstyle={inner sep=0pt}](6)(2);
			\Edge[ style = {->,> = latex',pos = 0.7,dashed,out=150,in=30},,color=black
			, labelstyle={inner sep=0pt}](2)(4);
			\Edge[ style = {->,> = latex',pos = 0.7,out=-30,in=-150},,color=black
			, labelstyle={inner sep=0pt}](4)(2);
			\Edge[ style = {->,> = latex',pos = 0.3,out=60,in=-60},,color=black
			, labelstyle={inner sep=0pt}](4)(3);
			\Edge[ style = {->,> = latex',pos = 0.3,out=-120,in=120},,color=black
			, labelstyle={inner sep=0pt}](3)(4);
			\Loop[style={> = latex',->, out=60, in=120,line width=0.8pt}, dist=1.5cm](1)
			\Loop[style={> = latex',->, out=-60, in=-120,line width=0.8pt,dashed}, dist=1.5cm](2)
			\Loop[style={> = latex',->, out=-120, in=180,line width=0.8pt}, dist=1.5cm](4)
			\end{tikzpicture}
			\caption{\centering Node $5$ colors $1$ and $6$ colors $2$.}
			\label{g:G2}
		\end{subfigure}
		\vspace{.3cm}
		\begin{subfigure}{0.3\textwidth}
			\centering
			\begin{tikzpicture}[scale=0.5]
		\tikzset{VertexStyle1/.style = {shape = circle,
	        color=black,
	        fill=white!96!black,
	        minimum size=0.5cm,
			text = black,
			inner sep = 2pt,
			outer sep = 1pt,
			minimum size = 0.55cm}
		}		
			\tikzset{VertexStyle2/.style = {shape = circle,
	        color=black,
	        fill=black!96!white,
	        minimum size=0.5cm,
			text = white,
			inner sep = 2pt,
			outer sep = 1pt,
			minimum size = 0.55cm}
}
			\node[VertexStyle2,draw](1) at (0,3) {$\bf 1$};
			\node[VertexStyle2,draw](2) at (0,0) {$\bf 2$};
			\node[VertexStyle2,draw](3) at (-3,3) {$\bf 3$};
			\node[VertexStyle1,draw](4) at (-3,0) {$\bf 4$};
			\node[VertexStyle1,draw](5) at (3,3) {$\bf 5$};
			\node[VertexStyle1,draw](6) at (3,0) {$\bf 6$};
			\Edge[ style = {->,> = latex',pos = 0.2},,color=black
			, labelstyle={inner sep=0pt}](5)(1);
			\Edge[style = {->,> = latex',pos = 0.7},color=black
			, labelstyle={inner sep=0pt}](1)(3);
			\Edge[style = {->,> = latex',pos = 0.7},color=black
			, labelstyle={inner sep=0pt}](6)(2);
			\Edge[ style = {->,> = latex',pos = 0.7,dashed,out=150,in=30},,color=black
			, labelstyle={inner sep=0pt}](2)(4);
			\Edge[ style = {->,> = latex',pos = 0.7,out=-30,in=-150},,color=black
			, labelstyle={inner sep=0pt}](4)(2);
			\Edge[ style = {->,> = latex',pos = 0.3,out=60,in=-60},,color=black
			, labelstyle={inner sep=0pt}](4)(3);
			\Edge[ style = {->,> = latex',pos = 0.3,out=-120,in=120},,color=black
			, labelstyle={inner sep=0pt}](3)(4);
				\Loop[style={> = latex',->, out=60, in=120,line width=0.8pt}, dist=1.5cm](1)
			\Loop[style={> = latex',->, out=-60, in=-120,line width=0.8pt,dashed}, dist=1.5cm](2)
			\Loop[style={> = latex',->, out=-120, in=180,line width=0.8pt}, dist=1.5cm](4)
			\end{tikzpicture}
			\caption{\centering Node $1$ colors $3$.}
			\label{g:G3}
		\end{subfigure}
		\vspace{.3cm}
		\begin{subfigure}{0.3\textwidth}
			\centering
			\begin{tikzpicture}[scale=0.5]
		\tikzset{VertexStyle1/.style = {shape = circle,
	        color=black,
	        fill=white!96!black,
	        minimum size=0.5cm,
			text = black,
			inner sep = 2pt,
			outer sep = 1pt,
			minimum size = 0.55cm}
		}		
			\tikzset{VertexStyle2/.style = {shape = circle,
	        color=black,
	        fill=black!96!white,
	        minimum size=0.5cm,
			text = white,
			inner sep = 2pt,
			outer sep = 1pt,
			minimum size = 0.55cm}
}
			\node[VertexStyle2,draw](1) at (0,3) {$\bf 1$};
			\node[VertexStyle2,draw](2) at (0,0) {$\bf 2$};
			\node[VertexStyle2,draw](3) at (-3,3) {$\bf 3$};
			\node[VertexStyle2,draw](4) at (-3,0) {$\bf 4$};
			\node[VertexStyle1,draw](5) at (3,3) {$\bf 5$};
			\node[VertexStyle1,draw](6) at (3,0) {$\bf 6$};
			\Edge[ style = {->,> = latex',pos = 0.2},,color=black
			, labelstyle={inner sep=0pt}](5)(1);
			\Edge[style = {->,> = latex',pos = 0.7},color=black
			, labelstyle={inner sep=0pt}](1)(3);
			\Edge[style = {->,> = latex',pos = 0.7},color=black
			, labelstyle={inner sep=0pt}](6)(2);
			\Edge[ style = {->,> = latex',pos = 0.7,dashed,out=150,in=30},,color=black
			, labelstyle={inner sep=0pt}](2)(4);
			\Edge[ style = {->,> = latex',pos = 0.7,out=-30,in=-150},,color=black
			, labelstyle={inner sep=0pt}](4)(2);
			\Edge[ style = {->,> = latex',pos = 0.3,out=60,in=-60},,color=black
			, labelstyle={inner sep=0pt}](4)(3);
			\Edge[ style = {->,> = latex',pos = 0.3,out=-120,in=120},,color=black
			, labelstyle={inner sep=0pt}](3)(4);
			\Loop[style={> = latex',->, out=60, in=120,line width=0.8pt}, dist=1.5cm](1)
			\Loop[style={> = latex',->, out=-60, in=-120,line width=0.8pt,dashed}, dist=1.5cm](2)
			\Loop[style={> = latex',->, out=-120, in=180,line width=0.8pt}, dist=1.5cm](4)
			\end{tikzpicture}
			\caption{\centering Node $3$ colors $4$.}
			\label{g:G4}
		\end{subfigure}
		\caption{Example of a colorable graph.}
		
		\label{g:G}
	\end{figure}
\end{ex}

The following theorem now provides a necessary and sufficient graph theoretic condition for a given pattern matrix to have full row rank.
 \begin{thm}
 	\label{t:rank}
Let $\mathcal{M} \in \{0,\ast,?\}^{p \times q}$ be a pattern matrix with $p \leq q$. Then, $\calM$ has full row rank if and only if $G(\calM)$ is colorable.
 \end{thm}

It is clear from the definition of the color change rule that colorability of a given graph can be checked in polynomial time.

Finally, based on the rank test in Theorem \ref{t:algebraic} and the result in Theorem \ref{t:rank}, the following necessary and sufficient graph theoretic condition for controllability of a given structured system is obtained.
\begin{thm}\label{t:main}
	Let $\mathcal{A} \in \{0,\ast,?\}^{n \times n}$ and $\mathcal{B} \in \{0,\ast,?\}^{n \times m}$ be pattern matrices. Also, 
	let $\bar{\mathcal{A}}$ be obtained from $\mathcal{A}$ by modifying the diagonal entries of $\mathcal{A}$ as follows:
		\begin{equation}\label{e:barofX2}
		\bar{\mathcal{A}}_{ii} := \begin{cases}
		\ast & \text{if } \mathcal{A}_{ii} = 0,  \\
		? & \text{otherwise}.
		\end{cases}	
		\end{equation}
Then, the structured system $(\mathcal{A},\mathcal{B})$ is controllable if and only if both $G(\bbm \calA & \calB \ebm)$ and $G(\bbm \bar{\calA}& \calB \ebm)$ are colorable.
\end{thm}
	
As an example, we study controllability of the electrical circuit discussed in Example~\ref{ex:electriccircuit}. 

\begin{ex}
According to Example~\ref{ex:same circuit}, the electrical circuit depicted in Figure~\ref{g:electriccircuit} can be modelled as a structured system of the form \eqref{e:gss} where the pattern matrices $\calA$ and $\calB$ are given by:
$$
	\calA = \begin{bmatrix}
\ast & 0 & \ast\\
0 & 0 & \ast\\
? & \ast  & \ast
\end{bmatrix} \mbox{ and } \calB = 
\begin{bmatrix}
\ast & 0\\
0 & \ast \\ 
?& 0
\end{bmatrix}.
$$
Then, we obtain 
$$
\bar{\calA} = \begin{bmatrix}
? & 0 & \ast\\
0 & \ast & \ast\\
? & \ast  & ?
\end{bmatrix}.
$$
The graphs $G(\bbm \calA & \calB \ebm)$ and $G(\bbm \bar{\calA} & \calB \ebm)$ are depicted in Figure~\ref{g:patterncircuit} and Figure~\ref{g:patterncircuit2}, respectively. Both graphs are colorable. Indeed, node $5$ colors $2$, node $2$ colors $3$, and finally $3$ colors $1$ in both graphs. Therefore, the system $(\calA,\calB)$ is controllable by Theorem~ \ref{t:main}.

\begin{figure}[h!]
		\centering
		\centering
		\begin{subfigure}{0.3\textwidth}
			\centering
		\begin{tikzpicture}[scale=0.7]
		\tikzset{VertexStyle1/.style = {shape = circle,
	        color=black,
	        fill=white!93!black,
	        minimum size=0.5cm,
			text = black,
			inner sep = 2pt,
			outer sep = 1pt,
			minimum size = 0.55cm}
		}	
			\node[VertexStyle1,draw](1) at (0,2) {$\bf 1$};
			\node[VertexStyle1,draw](3) at (0,-2) {$\bf 3$};
			\node[VertexStyle1,draw](2) at (2.5,0) {$\bf 2$};
			\node[VertexStyle1,draw](4) at (-2.5,0) {$\bf 4$};
			\node[VertexStyle1,draw](5) at (5,0) {$\bf 5$};
		
		\Edge[ style = {->,> = latex',pos = 0.2,dashed,, out=-120, in=120,},color=black, labelstyle={inner sep=0pt}](1)(3);
		\Loop[style={> = latex',->, out=60, in=120,line width=0.8pt}, dist=1.5cm](1)
		\Edge[style = {->,> = latex',pos = 0.7,out=-90, in=0,},color=black
		, labelstyle={inner sep=0pt}](2)(3);
		\Edge[ style = {->,> = latex',pos = 0.7},color=black
		, labelstyle={inner sep=0pt}](3)(1);
		\Edge[ style = {->,> = latex',pos = 0.3},,color=black
		, labelstyle={inner sep=0pt}](3)(2);
		\Loop[style={> = latex',->, out=-60, in=-120,line width=0.8pt}, dist=1.5cm](3);		
		\Edge[ style = {->,> = latex',pos = 0.3},,color=black
		, labelstyle={inner sep=0pt}](4)(1);
		\Edge[ style = {->,> = latex',pos = 0.3,dashed},,color=black
		, labelstyle={inner sep=0pt}](4)(3);
		\Edge[ style = {->,> = latex',pos = 0.3},color=black
		, labelstyle={inner sep=0pt} ](5)(2);
		\end{tikzpicture}
		\caption{ The graph $G(\bbm \calA & \calB \ebm)$.}
		\label{g:patterncircuit}
		\end{subfigure}
			\begin{subfigure}{0.3\textwidth}
			\centering
			\begin{tikzpicture}[scale=0.7]
		\tikzset{VertexStyle1/.style = {shape = circle,
	        color=black,
	        fill=white!93!black,
	        minimum size=0.5cm,
			text = black,
			inner sep = 2pt,
			outer sep = 1pt,
			minimum size = 0.55cm}
		}	
			\node[VertexStyle1,draw](1) at (0,2) {$\bf 1$};
			\node[VertexStyle1,draw](3) at (0,-2) {$\bf 3$};
			\node[VertexStyle1,draw](2) at (2.5,0) {$\bf 2$};
			\node[VertexStyle1,draw](4) at (-2.5,0) {$\bf 4$};
			\node[VertexStyle1,draw](5) at (5,0) {$\bf 5$};
			
			\Edge[ style = {->,> = latex',pos = 0.2,dashed,, out=-120, in=120,},color=black, labelstyle={inner sep=0pt}](1)(3);
			\Loop[style={> = latex',->, out=60, in=120,line width=0.8pt,dashed}, dist=1.5cm](1)
			\Loop[style={> = latex',->, out=60, in=120,line width=0.8pt}, dist=1.5cm](2)
			\Edge[style = {->,> = latex',pos = 0.7,out=-90, in=0,},color=black
			, labelstyle={inner sep=0pt}](2)(3);
			\Edge[ style = {->,> = latex',pos = 0.7},color=black
			, labelstyle={inner sep=0pt}](3)(1);
			\Edge[ style = {->,> = latex',pos = 0.3},,color=black
			, labelstyle={inner sep=0pt}](3)(2);
			\Loop[style={> = latex',->, out=-60, in=-120,line width=0.8pt,dashed}, dist=1.5cm](3);	
			\Edge[ style = {->,> = latex',pos = 0.3},,color=black
			, labelstyle={inner sep=0pt}](4)(1);
			\Edge[ style = {->,> = latex',pos = 0.3,dashed},,color=black
			, labelstyle={inner sep=0pt}](4)(3);
			\Edge[ style = {->,> = latex',pos = 0.3},color=black
			, labelstyle={inner sep=0pt} ](5)(2);
			\end{tikzpicture}
			\caption{ The graph $G(\bbm \bar{\calA} & \calB \ebm)$.}
			\label{g:patterncircuit2}
		\end{subfigure}
	\caption{ The  graphs associated with the circuit in Example~ \ref{ex:electriccircuit}.}
	\end{figure}
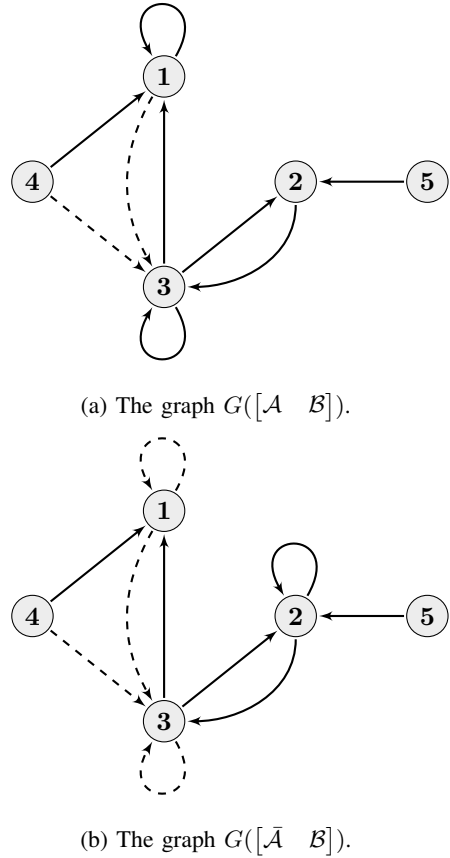	
\end{ex}

As a second example, we apply Theorem ~\ref{t:main} to verify the controllability of the networked system in Example \ref{ex:network}.
\begin{ex}
The networked system in Example \ref{ex:network} can be represented as a structured system of the form \eqref{e:gss}, where the pattern matrices $\calA$ and $\calB$ are given by:
	$$
	\calA = \begin{bmatrix}
	\ast & 0 & \ast\\
	? & \ast & 0\\
	0 & \ast  & 0
	\end{bmatrix} \mbox{ and } \calB = 
	\begin{bmatrix}
	0\\
	0  \\ 
	\ast
	\end{bmatrix}.
	$$
Clearly, 
	$$
	\bar{\calA} = \begin{bmatrix}
	? & 0 & \ast\\
	? & ? & 0\\
	0 & \ast  & \ast
	\end{bmatrix}.
	$$
	The graphs $G(\bbm \calA & \calB \ebm)$ and $G(\bbm \bar{\calA} & \calB \ebm)$ are depicted in Figure~\ref{g:patternnetwork} and Figure~\ref{g:patternnetwork2}, respectively. The graph in Figure~\ref{g:patternnetwork} is colorable. Indeed, node $4$ colors $3$, node $2$ colors $2$, and finally $3$ colors $1$. However, the graph in Figure~\ref{g:patternnetwork2} is not colorable. Therefore, the system $(\calA,\calB)$ is not controllable. However, if we would know that the edge $(1,2)$ does exist in the graph, i.e. if $\mathcal{A}_{21} = \ast$, then it can be verified that $(\mathcal{A},\mathcal{B})$ is controllable. 
	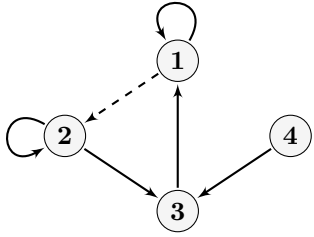
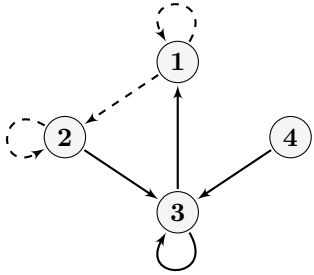
\begin{figure}[h!]
		\centering
		\centering
		\begin{subfigure}{0.3\textwidth}
			\centering
		\centering
		\begin{tikzpicture}[scale=0.5]
		\tikzset{VertexStyle1/.style = {shape = circle,
				color=black,
				fill=white!96!black,
				minimum size=0.5cm,
				text = black,
				inner sep = 2pt,
				outer sep = 1pt,
				minimum size = 0.55cm}
		}	
		\tikzset{VertexStyle2/.style = {shape = circle,
				ballcolor = black!80!yellow,
				text = white,
				inner sep = 2pt,
				outer sep = 0pt,
				minimum size = 10 pt}}
		\node[VertexStyle1,draw](1) at (0,2) {$\bf 1$};
		\node[VertexStyle1,draw](2) at (-3,0) {$\bf 2$};
		\node[VertexStyle1,draw](3) at (0,-2) {$\bf 3$};
		\node[VertexStyle1,draw] (4) at (3,0) {$\bf 4$};
		%
		\Edge[ style = {->,> = latex',pos = 0.2},color=black
		, labelstyle={inner sep=0pt}](2)(3);
		\Edge[ style = {->,> = latex',pos = 0.7, dashed},color=black, labelstyle={inner sep=0pt}](1)(2);
		\Loop[style={> = latex',->, out=150, in=-150,line width=0.8pt,color=black}, dist=1.5cm](2)
		\Edge[style = {->,> = latex',pos = 0.7},color=black
		, labelstyle={inner sep=0pt}](3)(1);	
		\Edge[ style = {->,> = latex',pos = 0.3},color=black
		, labelstyle={inner sep=0pt}](4)(3);
		\Loop[style={> = latex',->, out=60, in=120,line width=0.8pt,color=black}, dist=1.5cm](1);	
		\end{tikzpicture}
		\vspace*{3mm}
		\caption{The graph $G(\bbm \calA & \calB \ebm)$.}
		\label{g:patternnetwork}
		\end{subfigure}
		\begin{subfigure}{0.3\textwidth}
		\vspace*{3mm}
			\centering
			\begin{tikzpicture}[scale=0.5]
			\tikzset{VertexStyle1/.style = {shape = circle,
					color=black,
					fill=white!96!black,
					minimum size=0.5cm,
					text = black,
					inner sep = 2pt,
					outer sep = 1pt,
					minimum size = 0.55cm}
			}	
			\tikzset{VertexStyle2/.style = {shape = circle,
					ball color = black!80!yellow,
					text = white,
					inner sep = 2pt,
					outer sep = 0pt,
					minimum size = 10 pt}}
		\node[VertexStyle1,draw](1) at (0,2) {$\bf 1$};
		\node[VertexStyle1,draw](2) at (-3,0) {$\bf 2$};
		\node[VertexStyle1,draw](3) at (0,-2) {$\bf 3$};
		\node[VertexStyle1,draw] (4) at (3,0) {$\bf 4$};
			%
			\Edge[ style = {->,> = latex',pos = 0.2},color=black
			, labelstyle={inner sep=0pt}](2)(3);
			\Edge[ style = {->,> = latex',pos = 0.7, dashed},color=black, labelstyle={inner sep=0pt}](1)(2);
			\Loop[style={> = latex',->, out=150, in=-150,line width=0.8pt,color=black, dashed}, dist=1.5cm](2)
			\Edge[style = {->,> = latex',pos = 0.7},color=black
			, labelstyle={inner sep=0pt}](3)(1);	
			\Edge[ style = {->,> = latex',pos = 0.3},color=black
			, labelstyle={inner sep=0pt}](4)(3);
			\Loop[style={> = latex',->, out=60, in=120,line width=0.8pt,color=black,dashed}, dist=1.5cm](1);	
			\Loop[style={> = latex',->, out=-60, in=-120,line width=0.8pt,color=black}, dist=1.5cm](3);	
			\end{tikzpicture}
			\vspace*{2mm}
			\caption{ The graph $G(\bbm \bar{\calA} & \calB \ebm)$.}
			\label{g:patternnetwork2}
		\end{subfigure}
		\caption{The graphs associated with the network in Example~\ref{ex:network}.}
	\end{figure}
\end{ex}
 
By applying Theorem~\ref{t:main} to the special case discussed in Corollary~\ref{c:algebraic}, we obtain the following.

\begin{cor}\label{c:main}
Suppose that none of the diagonal entries of $\calA$ is zero. Let $\bar{\calA}$ be defined as in \eqref{e:barofX2}. Then, the system $(\mathcal{A},\mathcal{B})$ is controllable if and only if $G(\bbm \bar{\calA}& \calB \ebm)$ is colorable.
\end{cor}

To conclude this section, the results we have obtained for controllability lead to an interesting observation in the context of structural stabilizability. We say that a structured system $(\calA, \calB)$ is {\em stabilizable\/} if the pair $(A,B)$ is stabilizable for all $A \in \calP(\calA)$ and $B \in \calP(\calB)$.

For a single linear system, controllability implies stabilizability, whereas the reverse implication does not hold in general. Interestingly, for structured systems controllability and stabilizability do turn out to be equivalent, as stated next.

\begin{thm}\label{t:stabilizability}
The system $(\mathcal{A},\mathcal{B})$ is stabilizable if and only if it is controllable.
\end{thm}

\section{Discussion of existing results}\label{s:dis}
In this section, we compare our results with those existing in the literature. We begin with giving an account of the most relevant related work.

In the past, the strong structural controllability problem was studied almost exclusively (with the exception of \cite{MZC2014}) for systems of the form \eqref{e:gss} where the pattern matrices $\calA$ and $\calB$ do {\em not\/} contain $?$ entries, that is where $\mathcal{A} \in \{0,*\}^{n \times n}$ and $\mathcal{B} \in \{0,*\}^{n\times m}$. 
Within this line of research, the earliest work is \cite{MY1979} that considered the single-input case, i.e. $m=1$. The results of this paper were extended to the multi-input case in \cite{Bachmann1981}. The necessary and sufficient conditions (\cite[Thm. 1 and Thm. 2]{MY1979} and \cite[Satz 3]{Bachmann1981}) that these papers provide are graph theoretic in nature. 
For the same class of structured systems, but for the single input case, Olesky \emph{et al.} provided algebraic conditions for strong structural controllability\footnote{The authors use the terminology ``qualitative controllability" instead of ``strong structural controllability".} in \cite[Thm. 2.2, Thm. 2.4]{OMP1993}, which can also be interpreted in a graph theoretic context. Reinschke \emph{et al.} presented another graph theoretic test \cite[Thm. 1]{RSW1992} as well as an algebraic test \cite[Thm. 2]{RSW1992}. Later, Jarczyk \emph{et al.} pointed out that the graph theoretic test given in \cite{RSW1992} is erroneous \cite[Ex. 1]{JSA2011} and provided a correction \cite[Thm. 5]{JSA2011}.
The study of strong structural controllability has seen a recent revival in the context of networked systems. This line of research was initiated in \cite{CM2013} and followed up in the papers \cite{MZC2014} and \cite{TD2015}. These papers study also particular classes of systems of the form \eqref{e:gss}. More specifically,  \cite{CM2013} and \cite{TD2015} consider pattern matrices $\mathcal{A} \in \{0,*\}^{n \times n}$ and $\mathcal{B} \in \{0,*\}^{n\times m}$ with the additional assumption that $\calB$ is a pattern matrix with exactly one $*$ in each column and at most one $*$ in each row. The paper \cite{MZC2014} considers\footnote{In fact, \cite{MZC2014} considers only binary matrices $B$ in \eqref{e:gss}, that is $B\in\pset{0,1}^{n\times n}$, with exactly one $1$ in each column and at most one $1$ in each row. Since the image of $B$ would not be changed if $1$'s are replaced by $*$'s, considering binary matrices or $\pset{0,*}$- matrices with the same pattern do not make a difference in the study of controllability.} $\mathcal{A} \in \{0,*,?\}^{n \times n}$ and $\mathcal{B} \in \{0,*\}^{n\times m}$ with the additional assumption that all diagonal entries of $\calA$ are $?$, none of the off-diagonal entries is $?$, and $\calB$ is a pattern matrix with exactly one $*$ in each column and at most one $*$ in each row. The main results of the papers \cite{CM2013,MZC2014,TD2015} involve algebraic as well as  graph theoretic necessary and sufficient conditions for the classes they study. 
In the sequel, we will discuss how the results in the above-mentioned papers compare to the results in the present paper, in particular, with an eye towards algorithmic complexity as well as conceptual simplicity. 

\subsection{Graph theoretic conditions}

The graph theoretic conditions provided in \cite[Thm. 1]{MY1979} for the single-input case ($m=1$) and extended to the multi-input case in \cite[Satz 3]{Bachmann1981} are based on the graph $G=(V,E)$ associated with a pattern matrix $\bbm\calA & \calB\ebm$ where $\mathcal{A} \in \{0,*\}^{n \times n}$ and $\mathcal{B} \in \{0,*\}^{n\times m}$. 
Note that $V=\nset{n+m}$ in this case. The graph theoretic characterization in \cite[Satz 3]{Bachmann1981} (or in \cite[Thm. 1]{MY1979} if $m=1$) consists of three conditions. The first one requires checking the so-called accessibility of each node in $\nset{n}$ from the nodes in $\pset{n+1,n+2,\ldots,n+m}$. The remaining two conditions require checking certain relations for {\em all\/} subsets of $\nset{n}$. As such, the computational complexity of checking these conditions is {\em at least\/} exponential in $n$. 
Note that, in contrast, the computational complexity of checking the colorability conditions of our Theorem~\ref{t:main} is polynomial in $n$. 

The paper \cite{MY1979} provides another set of graph theoretic conditions, stated, more specifically, in \cite[Thm. 2]{MY1979} (only for the case $m=1$). As argued in \cite[p. 135]{MY1979}, this theorem performs better than \cite[Thm. 1]{MY1979} for sparse graphs. Essentially, the conditions given in \cite[Thm. 2]{MY1979} require checking the existence of a unique serial buds cactus as well as nonexistence of certain cycles within the graph $G$. How these conditions can be checked in an algorithmic manner is not clear, whereas the colorability conditions given in Theorem~\ref{t:main} can be checked by a simple algorithm.  

On top of the advantages of computational complexity, the conditions provided in Theorem~\ref{t:main} are more attractive because of their conceptual simplicity. Indeed, colorability is a simpler and more intuitive notion than those appearing in the results of  \cite{MY1979} and \cite{Bachmann1981}.

Yet another graph theoretical characterization is provided in \cite[Thm. 5]{JSA2011}. In order to verify the conditions of \cite[Thm. 5]{JSA2011}, one needs to check whether a unique spanning cycle family with certain properties exists in $n+m \choose n$ directed graphs obtained from the pattern matrices $\calA$ and $\calB$. Needless to say, checking the conditions of Theorem~\ref{t:main} is much easier than checking these conditions.

Also in the context of networked systems, graph theoretic conditions for strong structural controllability have been obtained (see e.g. \cite{CM2013,MZC2014,TD2015}).
To elaborate further on the relationship between the work on networked systems and our work, we first need to explain the framework of the papers \cite{CM2013,MZC2014,TD2015}. 
The starting point of these papers is a directed graph $H=(W,F)$ where $W=\nset{n}$ denotes the node set and $F$ the edge set. The graphs considered in \cite{CM2013,TD2015} are so-called loop graphs, that are graphs which are allowed to contain self-loops,  whereas \cite{MZC2014} does not allow self-loops. Apart from the graph $H$, these papers consider a subset of the node set $W$, the so-called leader set, say $W_L=\pset{w_1,w_2,\ldots,w_m}$. Based on the graph $H$ and $W_L$, \cite{CM2013,MZC2014,TD2015} introduce systems of the form \eqref{e:gss} 
where the pattern matrix $\calB$ is defined by
\beq\label{e:b for net}
\calB_{ij}=\begin{cases}
*&\text{ if }i=w_j\\
0&\text{ otherwise}
\end{cases}
\eeq
for $i\in\nset{n}$, $j\in\nset{m}$.  In \cite{CM2013} and \cite{TD2015} the pattern matrix $\calA$  is defined by
\beq\label{e:a for net1}
\calA_{ij}=\begin{cases}
*&\text{ if }(j,i)\in F\\
0&\text{ otherwise}
\end{cases}
\eeq
whereas in \cite{MZC2014} the pattern matrix $\calA$  is defined by
\beq\label{e:a for net2}
\calA_{ij}=\begin{cases}
*&\text{ if }(j,i)\in F\\
?&\text{ if }i=j\\ 
0&\text{ otherwise}
\end{cases}
\eeq
for $i,j\in\nset{n}$.

In \cite{CM2013}, the authors first define two bipartite graphs obtained from the pattern matrices $\calA$ and $\calB$. Then, they show in \cite[Thm. 5]{CM2013} that $(\calA,\calB)$ is  strongly structurally controllable if and only if there exist so-called constrained matchings with certain properties in these bipartite graphs. Later, in \cite[Thm. 5.4]{TD2015} an equivalence between the existence of constrained matchings and so-called zero forcing sets for loop graphs was established. To explain this in more detail, we need to introduce the notion of zero forcing that was originally studied in the context of minimal rank problems (see e.g. \cite{H2010}). 

Let $H=(W,F)$ be a directed loop graph and $S\subseteq W$. Color all nodes in $S$ black and the others white.

If a node $i$ (of any color) has exactly one white out-neighbor $j$, we change the color of $j$ to black and write $i\to j$. If all the nodes in $W$ can be colored black by repeated application of this color change rule, we say that $S$ is a {\em loopy zero forcing set} for $H$. Given a loopy zero forcing set, we can list the color changes in the order in which they were performed to color all nodes black. This list is called a chronological list of color changes.

In order to quote \cite[Thm. 5.5]{TD2015}, we need two more definitions. Define $W_{\mathrm{loop}} \subseteq W$ to be the subset of all nodes with self-loops and let $H^*$ be the graph obtained from $H$ by placing a self-loop at every node.

\begin{thm}\cite[Thm. 5.5]{TD2015} \label{t:TD} Let $H$ be a directed loop graph and $W_L$ be a leader set. Consider the pattern matrices defined in \eqref{e:b for net} and \eqref{e:a for net1}. Then, the structured system $(\mathcal{A},\mathcal{B})$ is controllable if and only if the following conditions hold:
	\begin{enumerate}[label = \arabic*.]
		\item $W_L$ is a loopy zero forcing set for $H$.
		\item $W_L$ is a loopy zero forcing set for $H^*$ for which there is a chronological list of color changes that does not contain a color change of the form $i\to i$ with $i\in W_{\mathrm{loop}}$.		
	\end{enumerate}
\end{thm}

A result similar to this theorem was obtained in \cite{MZC2014} for controllability of pattern matrices defined by \eqref{e:b for net} and \eqref{e:a for net2} that are obtained from a graph $H$ {\em without\/} self-loops. However, in order to deal with this class of pattern matrices, \cite{MZC2014} introduces a slightly different notion of zero forcing to be defined below.

Let $H=(W,F)$ be a directed graph without self-loops and $S\subseteq W$. Color all nodes in $S$ black and the others white.
If a {\em black\/} node $i$ has exactly one white out-neighbor $j$, we change the color of $j$ to black. If all the nodes in $W$ can be colored black by repeated application of this color change rule, we say that $S$ is a {\em ordinary zero forcing set} for $H$.

We now state the graph theoretic characterization of controllability established in \cite{MZC2014}. 
\begin{thm}\cite[Thm. IV.4]{MZC2014} \label{t:MZC} Let $H$ be a directed graph without self-loops and $W_L$ be a leader set. Consider the pattern matrices given by \eqref{e:b for net} and \eqref{e:a for net2}. Then, the structured system $(\mathcal{A},\mathcal{B})$ is controllable if and only if $W_L$ is an ordinary zero forcing set for $H$.
\end{thm}

Even though Theorems~\ref{t:TD} and \ref{t:MZC} present conditions that are similar in nature, it is not possible to compare these results immediately as they deal with two different and non-overlapping system classes. Indeed, the pattern matrices considered in \cite{TD2015} (given by \eqref{e:a for net1}) {\em do not contain\/} any $?$ entries whereas those studied in \cite{MZC2014} (given by \eqref{e:a for net2}) contain {\em only\/} $?$ entries on their diagonals.

Next, we will show that the conditions of Theorem~\ref{t:main} are equivalent to those of Theorems~\ref{t:TD} and \ref{t:MZC} if specialized to the corresponding pattern matrices. This will shed light on the relationship between these results based on the different zero forcing notions. 

We start with Theorem~\ref{t:TD}. According to our color change rule, the nodes belonging to $W_L$ will be colored black in both $G(\bbm \calA&\calB \ebm)$ and $G(\bbm\bar{\calA}&\calB\ebm)$ because $\calB$ is a pattern matrix with structure defined by \eqref{e:b for net}. Since $\calA$ does not contain $?$ entries, $G(\bbm \calA&\calB \ebm)$ is colorable if and only if $W_L$ is a loopy zero forcing set for $G(\calA)$.		
 By noting that $H=G(\calA)$, we see that the first condition in Theorem~\ref{t:main} is equivalent to that of Theorem~\ref{t:TD}. Now, let the pattern matrix $\calA^*$ be such that $H^*=G(\calA^*)$. Since $W_{\mathrm{loop}}=\set{i}{\bar{\calA}_{ii} = \mbox{?}}$, we see that $G(\bbm\bar{\calA}&\calB\ebm)$ is colorable if and only if the second condition of Theorem \ref{t:TD} holds. Thus, the second condition of Theorem~\ref{t:main} is equivalent to that of Theorem~\ref{t:TD}. 

Now, we turn attention to Theorem~\ref{t:MZC}. It follows from \eqref{e:barofX2} and \eqref{e:a for net2} that $\bar{\calA}=\calA$, i.e., graphs $G(\bbm \bar{\calA}&\calB \ebm)$ and $G(\bbm \calA &\calB \ebm)$ are the same. 
As in the discussion above, the nodes belonging to $W_L$ will be colored black in $G(\bbm \bar{\calA}&\calB \ebm)$ because $\calB$ is a pattern matrix with structure defined by \eqref{e:b for net}. According to our color change rule, a white node can never color any other white node in $G(\bbm \bar{\calA}&\calB \ebm)$ since $(i,i)\in E_?$ for every node $i$ of $G(\bar{\calA})$. This means that $G(\bbm \bar{\calA}&\calB \ebm)$ is colorable if and only if $W_L$ is an ordinary zero forcing set for $G(\bar{\calA})$. By noting that $H = G(\calA) =G(\bar{\calA})$, we see that the conditions in Theorem~\ref{t:main} are equivalent to the single condition of Theorem~\ref{t:MZC}.

\subsection{Algebraic conditions}

In this subsection, we will compare our rank tests for strong structural controllability with those provided in \cite{RSW1992,CM2013,MZC2014}. More precisely, we will show that the rank tests in Theorem \ref{t:algebraic} reduce to those in \cite{RSW1992,CM2013,MZC2014} for the corresponding special cases of pattern matrices.

An algebraic condition for controllability of $(\calA,\calB)$ was provided in \cite[Thm. 2]{RSW1992} for $\calA \in \{0,\ast\}^{n \times n}$ and $\calB \in \{0,\ast\}^{n \times m}$. 
Later, these conditions were reformulated in \cite[Thm. 3]{CM2013}. These conditions rely on a matrix property that will be defined next for pattern matrices that may also contain $?$ entries. 
\begin{defn} \label{d:reRSW1992}
	Consider a pattern matrix $\mathcal{M} \in \{0,\ast,?\}^{p \times q}$ with $p \leq q$. The matrix $\mathcal{M}$ is said to be of Form \Rmnum{3} if there exist two permutation matrices $P_1$ and $P_2$ such that
	\begin{equation}\label{e:reform3}
	P_1\mathcal{M} P_2 = 
	\begin{bmatrix}
	\otimes & \cdots & \otimes & \ast & 0 & \cdots & 0 \\
	\vdots &  & \vdots & \ddots & \ddots &  \ddots & \vdots\\
	\otimes & \cdots & \otimes & \cdots & \otimes & \ast & 0\\
	\otimes & \cdots & \otimes& \cdots & \otimes &\otimes & \ast
	\end{bmatrix},
	\end{equation}
where the symbol $\otimes$ indicates an entry that can be either $0$, $\ast$ or $?$.
\end{defn}

The above-mentioned algebraic conditions are stated next.
\begin{thm}\cite{CM2013}\label{t:CMalg}
Let $\mathcal{A} \in \{0,\ast\}^{n \times n}$ and $\mathcal{B} \in \{0,\ast\}^{n \times m}$ be two pattern matrices. Also, let $\calA^\ast$ be the pattern matrix obtained from $\calA$ by replacing all diagonal entries by $\ast$. The system $(\mathcal{A},\mathcal{B})$ is controllable if and only if the following two conditions hold:
	\begin{enumerate}[label = \arabic*.]
		\item The matrix $\bbm \calA & \calB \ebm$ is of Form \Rmnum{3}.\\[-3mm]
		\item The matrix $\bbm \calA^\ast &\!\!\! \calB \ebm$ is of Form \Rmnum{3} with the additional pro\-perty that $\ast$ entries appearing in \eqref{e:reform3} do not originate from diagonal elements in $\mathcal{A}$ that are $\ast$ entries.
	\end{enumerate}
\end{thm}

It can be shown that our algebraic conditions in Theorem \ref{t:algebraic} are equivalent to those in Theorem \ref{t:CMalg} for the special case of pattern matrices that only contain $0$ and $\ast$ entries. Recall that it follows from Theorem \ref{t:algebraic} that $(\mathcal{A},\mathcal{B})$ is controllable if and only if both $\bbm \calA & \calB \ebm$ and  $\bbm \bar{\calA} & \calB \ebm$ have full row rank, where $\bar{\calA}$ is given in \eqref{e:barofX2}. To relate our algebraic conditions with the ones in Theorem~\ref{t:CMalg}, we need the following lemma.

\begin{lem}\label{l:algebraic}
	Let $\mathcal{M} \in \{0,\ast,?\}^{p \times q}$ with $p \leq q$. Then, $\calM$ has full row rank if and only if $\mathcal{M}$ is of Form \Rmnum{3}.
\end{lem}

From Lemma \ref{l:algebraic} it immediately follows that $\bbm \calA & \calB \ebm$ has full row rank if and only if $\bbm \calA&\calB\ebm$ is of Form \Rmnum{3}. 
Hence, the first condition of Theorem \ref{t:algebraic} is equivalent to that of Theorem \ref{t:CMalg}. We will now also show that $\bbm \bar{\calA} & \calB \ebm$ has full row rank if and only if the second condition of Theorem \ref{t:CMalg} holds. 
From Lemma \ref{l:algebraic}, we have that $\bbm \bar{\calA} & \calB \ebm$ has full row rank if and only if $\bbm\mathcal{\bar{A}}&\mathcal{B}\ebm$ is of Form \Rmnum{3}.	
By definition of $\bar{\calA}$ and $\calA^\ast$, it follows that $\bar{\calA}_{ij} = \calA^\ast_{ij}$ for all $i \neq j$. If $\calA_{ii} = 0$ then both $\bar{\calA}_{ii} = \ast$ and $\calA^*_{ii} = *$. On the other hand, if $\calA_{ii} = *$ then $\bar{\calA}_{ii} = ?$ and $\calA^*_{ii} = *$. 
To sum up, $\bar{\calA}_{ij} \neq \calA^\ast_{ij}$ if and only if $i = j$ and $\calA_{ii} = \ast$.
In other words, all entries of $\bar{\calA}$ and $\calA^*$ are \emph{the same}, except for those that correspond to the diagonal elements of $\calA$ that are $\ast$ entries.   
Hence, there exist two permutation matrices $P_1$ and $P_2$ such that all entries of the matrices $P_1 \bbm \bar{\calA} & \calB \\ \ebm P_2$ and $ P_1 \bbm \calA^\ast & \calB \\ \ebm P_2$ are the same, except those that originate from diagonal elements of $\calA$ that are $\ast$ entries.
This implies that $\bbm \bar{\calA} & \calB \ebm$ is of Form \Rmnum{3} if and only if $\bbm \calA^\ast & \calB \ebm$ is of Form \Rmnum{3} with the additional property that the $\ast$ entries in \eqref{e:reform3} do not originate from diagonal elements in $\mathcal{A}$ that are $\ast$ entries. In other words, the second conditions of Theorem \ref{t:algebraic} and \ref{t:CMalg} are equivalent. 
Since also the first conditions in these theorems are equivalent, we conclude that the algebraic conditions in Theorem \ref{t:algebraic} are equivalent to those in Theorem \ref{t:CMalg} for the special case in which $\mathcal{A} \in \{0,\ast\}^{n \times n}$ and $\mathcal{B} \in \{0,\ast\}^{n \times m}$.

A different algebraic condition was introduced in \cite{MZC2014} for systems defined on simple directed graphs. The pattern matrices of such systems can be represented by $\calA$ and $\calB$ given by \eqref{e:a for net2} and \eqref{e:b for net}, respectively. 
The algebraic condition referred to above is then stated as follows.
\begin{thm}\cite[Lem. \Rmnum{4}.1]{MZC2014} \label{t:alMZC}
    	 Consider the pattern matrices $\calA$ and $\calB$ given by \eqref{e:a for net2} and \eqref{e:b for net}, respectively. Then,  $(\mathcal{A},\mathcal{B})$ is controllable if and only if $\bbm \calA & \calB\ebm$ has full row rank.
    \end{thm}
In order to see that this theorem follows from Corollary \ref{c:algebraic}, note that $\mathcal{A} = \bar{\mathcal{A}}$ since all diagonal entries of $\mathcal{A}$ are $?$'s.

	\section{Proofs}\label{S:proof}

\subsection{Proof of Theorem \ref{t:algebraic}}	
		To prove the `only if' part, assume that $(\mathcal{A},\mathcal{B})$ is controllable. By the Hautus test \cite[Thm. 3.13]{TSH2012} and the definition of strong structural controllability, it follows that $\begin{bmatrix} A- \lambda I & B
	\end{bmatrix}$ has full row rank for all $(A,B) \in \mathcal{P}(\mathcal{A}) \times \mathcal{P}(\mathcal{B})$ and all $\lambda \in \mathbb{C}$. By substitution of $\lambda = 0$ we conclude that condition 1 is satisfied. To prove that condition 2 also holds, suppose that $x^T \begin{bmatrix}
	\bar{A} &  B
	\end{bmatrix} = 0$ for some pair $(\bar{A},B) \in \mathcal{P}(\mathcal{\bar{A}}) \times \mathcal{P}(\mathcal{B})$ and $x \in \mathbb{R}^n$. We want to prove that $x = 0$. Let $\alpha \in \mathbb{R}$ be a nonzero real number such that  
	\begin{equation*}
	\alpha \not\in \{\bar{A}_{ii} \mid i \text{ is such that } \mathcal{A}_{ii} = \ast \}.
	\end{equation*}
	Then, define a nonsingular diagonal matrix $X \in \mathbb{R}^{n \times n}$ as
	\begin{equation*}
	X_{ii} = \begin{cases}
	1 & \text{if } \bar{\mathcal{A}}_{ii} = \: ? \\
	\alpha \slash \bar{A}_{ii} & \text{if } \bar{\mathcal{A}}_{ii} = \ast.
	\end{cases}
	\end{equation*}
	It is clear that $\bar{A}X \in \mathcal{P}({\bar{\mathcal{A}}})$ and $x^T \begin{bmatrix}
	\bar{A}X & B
	\end{bmatrix} = 0$. Furthermore, by the choice of $\alpha$ and $X$ we obtain $\hat{A} := \bar{A}X - \alpha I \in \mathcal{P}({\mathcal{A}})$. By assumption, $\begin{bmatrix} \hat{A} + \alpha I & B \end{bmatrix}$ has full row rank (by substitution of $\lambda = - \alpha$). In other words, $\begin{bmatrix}
	\bar{A}X & B
	\end{bmatrix}$ has full row rank and therefore $x = 0$. We conclude that condition 2 is satisfied.
	
	To prove the `if' part, assume that conditions 1 and 2 are satisfied. Suppose that 
$$
	z^H\begin{bmatrix}
	A - \lambda I & B
	\end{bmatrix} = 0
$$
	for some $(A,B) \in \mathcal{P}(\mathcal{A}) \times \mathcal{P}(\mathcal{B})$ and $(\lambda,z) \in \mathbb{C} \times \mathbb{C}^n$, and $z^H$ denotes the conjugate transpose of $z$. We want to prove that $z = 0$. Note that if $\lambda = 0$, it readily follows that $z = 0$ by condition 1.  Therefore, it remains to be shown that $z = 0$ if $ \lambda \neq 0$. To this end, write $z = \xi + j\eta$, where $\xi,\eta \in \mathbb{R}^n$ and $j$ denotes the imaginary unit. 
Next, let $\alpha \in \mathbb{R}$ be a nonzero real number such that
	\begin{equation*}
	\alpha \not\in \left \{ -\frac{\xi_i}{\eta_i} \mid \eta_i \neq 0 \right \} \cup \left \{ -\frac{(\xi^T A)_i}{(\eta^T A)_i} \mid (\eta^T A)_i \neq 0 \right \}.
	\end{equation*}
	We define $x := \xi+\alpha \eta$. Now, we claim that
	\begin{enumerate}[label=(\alph*)]
		\item  $x_i = 0$ if and only if $z_i = 0$. 
		\item $x_i = 0$ if and only if $(x^T A)_i = 0$. 
	\end{enumerate}
	Note that (a) follows directly from the definition of $x$ and the choice of $\alpha$. To prove the `only if' part of (b), suppose that $x_i = 0$. By (a), this implies that $z_i = 0$. Since $z^H A = \lambda z^H$, we see that $(z^H A)_i = 0$. Equivalently, $((\xi^T - j\eta^T)A)_i = 0$. Therefore, both $(\xi^T A)_i = 0$ and $(\eta^T B)_i = 0$. We conclude that $(x^T A)_i = ((\xi^T + \alpha \eta^T)A)_i = 0$. 
	
	To prove the `if' part of (b), suppose that $(x^T A)_i = 0$. This means that $((\xi^T + \alpha \eta^T) A)_i = 0$. Equivalently, $(\xi^T A)_i + \alpha (\eta^T A)_i = 0$. By the choice of $\alpha$, this implies that $(\xi^T A)_i = (\eta^T A)_i = 0$. We conclude that $(z^H A
	)_i = 0$. Recall that $z^H A = \lambda z^H$, where $\lambda$ was assumed to be \emph{nonzero}. This implies that $z_i = 0$. Again, using (a) we conclude that $x_i = 0$. This proves (b).
	
	Next, we define the diagonal matrix $X' \in \mathbb{R}^{n \times n}$ as 
	\begin{equation*}
	X'_{ii} = \begin{cases}
	1 & \text{if } x_i = 0 \\
	\frac{(x^TA)_i}{x_i} & \text{otherwise}.
	\end{cases}
	\end{equation*}
	We know that $X'$ is nonsingular by (b). By definition of $X'$ we have $x^T A = x^T X'$. Furthermore, as $z^H B = 0$ we obtain $\xi^T B = \eta^T B = 0$ and therefore $x^T B = 0$. Hence $x^T \begin{bmatrix}
	A-X' & B
	\end{bmatrix} = 0$. Since $X'$ is nonsingular, $A - X' \in \mathcal{P}({\bar{\mathcal{A}}})$. By condition 2, this means that $x = 0$. Finally, we conclude that $z = 0$ using (a). 
	\EP
\subsection{Proof of Theorem \ref{t:rank}}

To prove Theorem \ref{t:rank}, we need the following auxiliary result.

\blem \label{l:fullrank}
Let $\mathcal{M} \in \{0,\ast,?\}^{p \times q}$ be a pattern matrix with $p \leq q$. Consider the directed graph $G(\calM)$. Suppose that each node is colored white or black. Let $D\in \mathbb{R}^{p \times p}$ be the diagonal matrix defined by
$$
D_{kk} = \begin{cases} 1 & \text{if node } k \text{ is black}, \\
0 & \text{otherwise.}
\end{cases}
$$
Suppose further that $j\in\nset{p}$ is a node for which there exists a node $i\in\nset{p}$, possibly identical to $j$, such that $j$ is the only white out-neighbor of $i$ and $(i,j)\in E_*$. Then for all $M \in \mathcal{P}(\mathcal{M})$ we have that $\begin{bmatrix} M & D \end{bmatrix}$ has full row rank if and only if $\begin{bmatrix} M & D + e_je_j^T \end{bmatrix}$ has full row rank where $e_j$ denotes the $j$th column of $I$.
\elem
\begin{IEEEproof}
	The `only if' part is trivial. To prove the `if' part, suppose that $M \in \mathcal{P}(\mathcal{M})$ and $\begin{bmatrix}
	M & D + e_je_j^T \end{bmatrix}$ has full row rank. 
	Let $z\in\R^p$ be such that $z^T\bbm M & D\ebm=0$. Our aim is to show that $z_j=0$. 
	Indeed, if $z_j$ is zero then $z^T\bbm M & D+ e_je_j^T \ebm=z^T\bbm M & D\ebm=0$ and hence $z$ must be zero. This would prove that $\bbm M & D\ebm$ has full row rank. 
   We  will distinguish two cases: $i = j$ and $i \neq j$.
Suppose first that $i = j$. Let $\beta,\omega\subseteq\nset{p}$  be defined as the index sets $\beta=\set{k}{k\neq j\text{ and } k \text{ is black}}$ and $\omega=\set{\ell}{\ell\neq j\text{ and } \ell \text{ is white}}$. 
	In the sequel, to simplify the notations, for a given vector $z \in\R^p$ and a given index set $\alpha\subseteq\nset{p}$,  we define $z_{\alpha} := \{x \in \R^{|\alpha|} \mid x_{i} = z_{\alpha(i)}, i \in  \nset{|\alpha|}\}$, where $|\alpha|$ is the cardinality of $\alpha$.
From $z^TM=0$, we get
\beq\label{e:index sets and such}
z_jM_{jj}+z_{\beta}^TM_{\beta j}+z_{\omega}^TM_{\omega j}=0.
\eeq
Since $j$ is the only white out-neighbor of itself, we must have that $M_{jj}$ is nonzero and that $M_{\omega j}$ is a zero vector. Moreover, it follows from $z^TD=0$ that $z_\beta$ must a zero vector. Therefore, \eqref{e:index sets and such} implies that $z_j$ must be zero.

Next, suppose that $i \neq j$. Let $\beta,\omega\subseteq\nset{p}$  be defined as the index sets $\beta=\set{k}{k\neq i,\,k\neq j,\text{ and } k \text{ is black}}$ and $\omega=\set{\ell}{\ell\neq i,\,\ell\neq j,\text{ and } \ell \text{ is white}}$. From $z^TM=0$, we now get
\beq\label{e:index sets and such-2}
z_iM_{ii}+z_jM_{ji}+z_{\beta}^TM_{\beta i}+z_{\omega}^TM_{\omega i}=0.
\eeq
Since $j$ is the only white out-neighbor of $i$, we must have that $M_{ji}$ is nonzero and that $M_{\omega i}$ is a zero vector. Moreover, it follows from $z^TD=0$ that $z_\beta$ must a zero vector. Therefore, \eqref{e:index sets and such-2} implies that
\beq\label{e:index sets and such-3}
z_iM_{ii}+z_jM_{ji}=0.
\eeq	
Now, we distinguish two cases: $i$ is black and $i$ is white. If $i$ is black, then we have that $z_i$ is zero because $z^TD=0$. Therefore, \eqref{e:index sets and such-3} implies that $z_j=0$ as desired. Finally, if $i$ is white, then we have that $M_{ii}=0$ for otherwise $i$ would have two white out-neighbors. Again, \eqref{e:index sets and such-3} implies that $z_j$ is zero. This completes the proof.
\end{IEEEproof}
Now, we can give the proof of Theorem \ref{t:rank}.
\begin{IEEEproof}[Proof of Theorem \ref{t:rank}]
	To prove the `if' part, suppose that $G(\calM)$ is colorable. Let $M \in \mathcal{P}(\mathcal{M})$ . By repeated application of Lemma \ref{l:fullrank}, it follows that $M$ has full row rank if and only if $\begin{bmatrix}
	M & I
	\end{bmatrix}$ has full row rank, which is obviously true. Therefore, we conclude that $M$ has full row rank. 
	
	To prove the `only if' part, suppose that $\calM$ has full row rank but $G(\calM)$ is not colorable. Let $C$ be the set of nodes that are colored black by repeated application of the color change rule until no more color changes are possible. Then, $C$ is a {\em strict\/} subset of $\{1,2,\ldots,p\}$. Thus, possibly after reordering the nodes, we can partition $\mathcal{M}$ as
$$
	\mathcal{M} = \begin{bmatrix}
	\mathcal{M}_1 \\  \mathcal{M}_2
	\end{bmatrix},
$$
	where the rows of the matrix $ \mathcal{M}_1 $ correspond to the nodes in $C$ and the matrix $ \mathcal{M}_1$ correspond to the remaining white nodes. 		
	Note that $C = \emptyset$ means
	that $\mathcal{M}_2 = \mathcal{M}$ and $\mathcal{M}_1$ is absent. Since no more color changes are possible, there is no column of $\calM_2$ that has exactly one $\ast$ entry while all other entries are $0$. Therefore, for any column of $\mathcal{M}_2$, we have one of the following three cases:
	\begin{enumerate}[label = \alph*.]			
		\item All entries are $0$.
		\item There exists exactly one $?$ entry while all other entries are $0$.
		\item At least two entries belong to the set $\{\ast, ?\}$.
	\end{enumerate}
Consequently, there exists a matrix $M_2\in\calP(\calM_2)$ such that its column sums are zero, that is $\mathbbm{1}^TM_2=0$, where $\mathbbm{1}$ denotes the vector of ones of appropriate size. Take any $M_1\in\calP(\calM_1)$. Then
$$M = \begin{bmatrix}
	M_1 \\ M_2
	\end{bmatrix} \in \calP(\begin{bmatrix}
	\mathcal{M}_1 \\  \mathcal{M}_2
	\end{bmatrix})=\mathcal{P}(\mathcal{M})$$
satisfies
	\begin{equation*}
	\begin{bmatrix}
	0^T & \mathbbm{1}^T
	\end{bmatrix}
	\begin{bmatrix}
	M_1 \\ M_2
	\end{bmatrix}
	= 0.
	\end{equation*}
Hence, $M$ does not have full row rank and we have reached a contradiction. 
\end{IEEEproof}

\subsection{Proof of  Theorem~\ref{t:main}}
By Theorem~\ref{t:algebraic} and Theorem~\ref{t:rank}, we have that $\bbm \mathcal{A} & \mathcal{B}\ebm$ is controllable if and only if if and only if $G(\bbm \calA & \calB \ebm)$ and $G(\bbm \bar{\calA}& \calB \ebm)$ are colorable.
\EP

\subsection{Proof of Theorem~\ref{t:stabilizability}}\label{a:stabilizability}
The `if' part is evident. Therefore, it is enough to prove the `only if' part. Suppose that the system $(\calA,\calB)$ is stabilizable. Let $(A,B)\in\calP(\calA)\times\calP(\calB)$. Then, $(A,B)$ is stabilizable. Note that $A\in\calP(\calA)$ if and only if $-A\in\calP(\calA)$. Therefore, we have both $(A,B)$ and $(-A,B)$ stabilizable. It follows from the Hautus test for stabilizability (see e.g. \cite[Thm. 3.32]{TSH2012}) that $(A,B)$ is controllable. Consequently, the system $(\calA,\calB)$ is controllable.\EP
\subsection{Proof of Lemma \ref{l:algebraic}}\label{a:s2}
Since the `if' part is evident, it remains to prove the `only if' part. Suppose that  $\calM$ has full row rank. From Theorem \ref{t:rank}, it follows that  $G(\mathcal{M})$ is colorable. In particular, there exist $i\in\nset{q}$ and $j\in\nset{p}$ such that $\calM_{ji}=\ast$ and $\calM_{ki}=0$ for all $k\neq j$. Therefore, we can find permutation matrices $P_1'$ and $P_2'$ such that
$$
P_1'\calM P_2'=
\left[\begin{array}{c|c}
\calM' & \begin{matrix}0 \\ \vdots \\ 0\end{matrix}\\
\hline
\begin{matrix}\otimes&\cdots&\otimes\end{matrix}& \ast
\end{array}\right]
$$
where the symbol $\otimes$ indicates an entry that can be either $0$, $\ast$ or $?$. Note that $M$ has full row rank for all $M\in\calP(\calM)$ if and only if $M'$ has full row rank for all $M\in\calP(\calM')$. Therefore, repeated application of the argument above results in permutation matrices $P_1$ and $P_2$ such that
$$
	P_1\mathcal{M} P_2 = 
	\begin{bmatrix}
	\otimes & \cdots & \otimes & \ast & 0 & \cdots & 0 \\
	\vdots &  & \vdots & \ddots & \ddots &  \ddots & \vdots\\
	\otimes & \cdots & \otimes & \cdots & \otimes & \ast & 0\\
	\otimes & \cdots & \otimes& \cdots & \otimes &\otimes & \ast
	\end{bmatrix}.
$$
\EP

\section{Conclusions}\label{s:con}

In most of the existing literature on strong structural controllability of structured systems, a zero/nonzero structure of the system matrices is assumed to be given. However, in many physical systems or linear networked systems, apart from fixed zero entries and nonzero entries we need to allow a third kind of entries, namely those that can take arbitrary (zero or nonzero) values. To deal with this, we have extended the notion of zero/nonzero structure to what we have called zero/nonzero/arbitrary structure. 
We have formalized this more general class of structured systems using pattern matrices containing fixed zero, arbitrary nonzero and arbitrary entries.  
In this setup, we have established necessary and sufficient algebraic conditions for strong structural controllability of these systems in terms of full rank tests on two associated pattern matrices.
Moreover, a necessary and sufficient graph theoretic condition for a given pattern matrix to have full row rank has been given in terms of a new color change rule. 
We have then established a graph theoretic test for strong structural controllability of the new class of structured systems.
Finally, we have shown how our results generalize previous work. We have also shown that some existing results  \cite{TD2015,MZC2014} that are seemingly incomparable to ours, can be put in our framework, thus unveiling an overarching theory.

In addition to strong structural controllability, weak structural controllability and strong structural stabilizability of structured systems with zero/nonzero/arbitrary structures have been briefly analyzed. We have shown that weak structural controllability of our structured systems can be checked using tests that already exist in the literature. We have also shown that a structured system with zero/nonzero/arbitrary structure is strongly structurally stabilizable if and only if it is strongly structurally controllable. 

It would be interesting to adopt our new framework of structured systems to other problem areas in systems and control, such as network identification \cite{WTC2018} or fault detection and isolation \cite{REC2015}. This is left as a possibility for future research.

\ifCLASSOPTIONcaptionsoff
  \newpage
\fi

\bibliographystyle{IEEEtran}
\bibliography{myieee2019}

\vfill
\end{document}

%% file: ka-newcommands.tex
\let\leq\leqslant

\let\emptyset\varnothing

\newcommand{\calA}{\ensuremath{\mathcal{A}}}
\newcommand{\calB}{\ensuremath{\mathcal{B}}}

\newcommand{\calM}{\ensuremath{\mathcal{M}}}

\newcommand{\calP}{\ensuremath{\mathcal{P}}}

\newcommand{\bmat}{\begin{matrix}}
\newcommand{\emat}{\end{matrix}}
\newcommand{\bbm}{\begin{bmatrix}}
\newcommand{\ebm}{\end{bmatrix}}
\newcommand{\bbma}{\begin{bmatrix*}[r]}
\newcommand{\ebma}{\end{bmatrix*}}
\newcommand{\bpm}{\begin{pmatrix}}
\newcommand{\epm}{\end{pmatrix}}
\newcommand{\bvm}{\begin{vmatrix}}
\newcommand{\evm}{\end{vmatrix}}
\newcommand{\bse}{\begin{subequations}}
\newcommand{\ese}{\end{subequations}}
\newcommand{\beq}{\begin{equation}}
\newcommand{\eeq}{\end{equation}}
\newcommand{\ben}{\renewcommand{\labelenumi}{\arabic{enumi}.}
\renewcommand{\theenumi}{\arabic{enumi}}\begin{enumerate}}

\newcommand{\een}{\end{enumerate}}

\newcommand{\beni}{\renewcommand{\labelenumi}{\roman{enumi}.}
\renewcommand{\theenumi}{\roman{enumi}}\begin{enumerate}}

\newcommand{\eeni}{\end{enumerate}}

\newcommand{\bena}{\renewcommand{\labelenumi}{\alph{enumi}.}
\renewcommand{\theenumi}{\alph{enumi}}\begin{enumerate}}

\newcommand{\eena}{\end{enumerate}}

\newcommand{\bit}{\begin{itemize}}
\newcommand{\eit}{\end{itemize}}
\newcommand{\bthe}{\begin{theorem}}
\newcommand{\ethe}{\end{theorem}}
\newcommand{\blem}{\begin{lemma}}
\newcommand{\elem}{\end{lemma}}
\newcommand{\bprop}{\begin{proposition}}
\newcommand{\eprop}{\end{proposition}}
\newcommand{\bex}{\begin{example}}
\newcommand{\eex}{\end{example}}
\newcommand{\bas}{\begin{assumption}}
\newcommand{\eas}{\end{assumption}}
\newcommand{\bre}{\begin{remark}}
\newcommand{\ere}{\end{remark}}
\newcommand{\bcor}{\begin{corollary}}
\newcommand{\ecor}{\end{corollary}}
\newcommand{\bdfn}{\begin{definition}}
\newcommand{\edfn}{\end{definition}}
\newcommand{\bcon}{\begin{conjecture}}
\newcommand{\econ}{\end{conjecture}}

\newcommand{\pset}[1]{\ensuremath{\{#1\}}}
\newcommand{\nset}[1]{\ensuremath{\{1,2,\ldots,#1\}}}

\newcommand{\set}[2]{\ensuremath{\{#1\mid #2\}}}

\newcommand{\R}{\ensuremath{\mathbb R}}

\newcommand{\EP}{\hspace*{\fill} $\blacksquare$\bigskip\noindent}